\documentclass[USenglish]{article}	

\usepackage[utf8]{inputenc}				
\usepackage[big,online]{dgruyter}	
\usepackage{lmodern} 
\usepackage{microtype}
\usepackage[numbers,square,sort&compress]{natbib}

\theoremstyle{dgthm}
\newtheorem{theorem}{Theorem}
\newtheorem{corollary}{Corollary}
\newtheorem{proposition}{Proposition}

\theoremstyle{dgdef}

\newtheorem{remark}{Remark}

\newcommand\slo{s}
\newcommand\sslo{\check{\slo}}
\newcommand\slowtil{\widetilde{\slo}}
\newcommand\slotil{\partial_z\slo}
\newcommand\ttf{\mathfrak{t}}

\newcommand\nlc{\eta}
\newcommand\nlcc{\check{\eta}}
\newcommand\Xslo{X_\slo}
\newcommand\Xnlc{X_\nlc}
\newcommand\dslo{\underline{\delta\slo}}
\newcommand\dnlc{\underline{\delta\nlc}}
\newcommand\dphi{\underline{\delta\phi}}
\newcommand\dpsi{\underline{\delta\psi}}
\newcommand\dq{\underline{\delta\vec{q}}}
\newcommand\dsslo{\underline{\delta\sslo}}
\newcommand\dnlcc{\underline{\delta\nlcc}}
\newcommand\dpsic{\underline{\delta\psic}}
\newcommand\phih{\hat{\phi}}
\newcommand\psih{\hat{\psi}}
\newcommand\psic{\check{\psi}}
\newcommand\pmeas{p^{\textup{meas}}}
\newcommand\zbot{y}
\def\Ltwotil{L^2(\tilde{\Omega}_\zbot)}
\def\Ltwotilbdy{L^2(\partial\tilde{\Omega}_\zbot)}
\newcommand\omegad{\omega_d}
\newcommand\omegap{\omega_p^3}
\renewcommand\kappa{\tilde{\omega}}

\newcommand{\uproman}[1]{\uppercase\expandafter{\romannumeral#1}}

\usepackage{color}
\begin{document}

	\articletype{Research Article}
	\received{Month	DD, YYYY}
	\revised{Month	DD, YYYY}
  \accepted{Month	DD, YYYY}
  \journalname{De~Gruyter~Journal}
  \journalyear{YYYY}
  \journalvolume{XX}
  \journalissue{X}
  \startpage{1}
  \aop
  \DOI{10.1515/sample-YYYY-XXXX}

\title{Simultaneous reconstruction of sound speed and nonlinearity parameter in a paraxial model of vibro-acoustography in frequency domain}
\runningtitle{On the inverse problem of vibro-acoustography}

\author*[1]{Barbara Kaltenbacher}
\author[2]{Teresa Rauscher}
\runningauthor{B.~Kaltenbacher and T.~Rauscher}
\affil[1]{\protect\raggedright 
University of Klagenfurt, Department of Mathematics, Klagenfurt, Austria, e-mail: barbara.kaltenbacher@aau.at}
\affil[2]{\protect\raggedright 
University of Klagenfurt, Department of Mathematics, Klagenfurt, Austria, e-mail: teresa.rauscher@aau.at}
	
	
\abstract{In this paper we consider the inverse problem of vibro-acoustography, a technique for enhancing ultrasound imaging by making use of nonlinear effects. It  amounts to determining two spatially variable coefficients in a system of PDEs describing propagation of two directed sound beams and the wave resulting from their nonlinear interaction.
To justify the use of Newton's method for solving this inverse problem, on one hand we verify well-definedeness and differentiability of the forward operator corresponding to two versions of the PDE model; on the other hand we consider an all-at-once formulation of the inverse problem and prove convergence of Newton's method for its solution.
}



\maketitle
	
\section{Introduction}
Recently, several approaches for enhancing ultrasound by means of nonlinear effects have been proposed. In this paper we consider vibro-acoustography which has originally been proposed in \cite{FatemiGreenleaf1998,FatemiGreenleaf1999} to achieve
the enhanced resolution by high frequency waves while avoiding the drawbacks of scattering from small inclusions and of stronger attenuation at higher frequencies.
The experiment for image acquisition 
basically consists of three parts: 
(i) Two ultrasound beams of high and slightly different frequencies $\omega_1$ and $\omega_2$ are excited at two parts $\Sigma_1$, $\Sigma_2$ of an array of piezoelectric  transducers (emitters). (ii) They interact nonlinearly in a focus region, thus exciting a wave that basically propagates at the difference frequency $\omega_1-\omega_2$ (which allows to avoid strong scattering and attenuation). (iii) The latter is eventually measured by a receiver array $\Gamma$.
We refer to, e.g., \cite{vibroacoustics,Part_linIP} for a graphical illustration.\\
The fact that the value of the nonlinearity parameter $\frac{B}{A}$ governing the interaction depends on the tissue properties and thus varies in space $\frac{B}{A}=\frac{B}{A}(x)$ yields a means of imaging. Likewise, also the sound speed typically exhibits spatial variation $c=c(x)$.
A modeling and simulation framework for this imaging technology has been devised in 
\cite{Malcolmetal2007,Malcolmetal2008}.

The aim of this paper is to study the inverse problem of identifying both $\frac{B}{A}(x)$ and $c(x)$ in an appropriate PDE model, in which these two quantities appear as coefficients. Some preliminary results on this parameter identification problem have been obtained in \cite{vibroacoustics}, however, without modeling the excitation waves as sound beams. The latter (along with an analysis of this model with variable coefficients) is one of the main novel contributions of this paper. Another one is a  convergence analysis of a Newton type method for solving the inverse problem. 

\medskip

The model we will use here can be derived from a general wave equation 
for an acoustic velocity potential $\phi$ given by
\begin{equation} \label{eq:der15}
	\partial_t^2 \phi - c^2 \Delta \phi = \partial_t |\nabla \phi |^2 - \frac{1}{2} (\nabla \phi) \cdot \nabla |\nabla \phi|^2 + 
\frac{B}{A}
\partial_t \phi \Delta \phi - 
\frac{B}{2A}
|\nabla \phi |^2 \Delta \phi,
\end{equation}
by means of the asymptotic ansatz
\begin{equation} \label{eq:dervib0}
	\tilde{\phi}(t,x,\varepsilon) = \varepsilon(\tilde{\phi}_1(t,x)+\tilde{\phi}_2(t,x)) + \varepsilon^2 \tilde{\psi}(t,x),
\end{equation}
with small $\varepsilon$ and 
a paraxial approximation (to take into account the fact that propagation is concentrated to a certain given axis) as well as transformation to frequency domain.
In \eqref{eq:der15}, $c$ is the speed of sound and 
$\frac{B}{A}$ the nonlinearity parameter.

\subsection*{Paraxial approximation}
We assume that the direction of propagation is the $x_1$-axis. Then, for the original  variables $(t,x)=(t, x_1, x')$, where $x'=(x_2, \ldots, x_d)$, {(and marking functions of these variables by a check)} we perform the paraxial change of variables 
{
\begin{equation}\label{parax}
\begin{aligned}
&	(\tau, z, y)=\left( t-\check{\ttf}, \tilde{\varepsilon} x_1, \sqrt{\tilde{\varepsilon}}x'\right), \quad
\textup{ with }\check{\ttf}=\check{\ttf}(x_1,x'), \quad \ttf=\ttf(z,\zbot) \\
&\textup{ such that }\frac{1}{c^2}=|\nabla_x\check{\ttf}|^2
=|\partial_{x_1}\check{\ttf}|^2+|\nabla_{x'}\check{\ttf}|^2
=\tilde{\varepsilon}^2|\partial_{z}\ttf|^2+\tilde{\varepsilon}|\nabla_{\zbot}\ttf|^2, 
\quad \check{\ttf}(0,x')=0, \ x'\in \tilde{\Omega}_{\zbot}
\end{aligned}
\end{equation}
where $\tilde{\varepsilon} \ll 1$, so $\tilde{\varepsilon} \ll \sqrt{\tilde{\varepsilon}}$, and in general $\tilde{\varepsilon}$ can be different from $\varepsilon$ that was introduced for the second order approximation in \eqref{eq:dervib0} but for the model to make sense the relation $\varepsilon \leq \tilde{\varepsilon} \leq \sqrt{\varepsilon}$ should hold cf. \cite{Part_linIP}.\footnote{Note that $\check{\ttf}(x)$ is the distance of $x$ from the boundary $\{0\}\times\tilde{\Omega}_{\zbot}$ in the Riemannian metric determined by $\frac{1}{c^2}$, cf., e.g. \cite[Section 7.2.3]{Evans:1998}}   
For the derivatives, this change of variables yields
\begin{equation}\label{parax_deriv}
\begin{aligned}
&	\partial_t=\partial_{\tau}, \quad 
\nabla_x=
\left( \tilde{\varepsilon}\bigl(-\partial_{z}\ttf \,\partial_{\tau} +  \partial_z\bigr), 
\sqrt{\tilde{\varepsilon}}\bigl(-\nabla_{\zbot}\ttf\, \partial_{\tau} + \nabla_{\zbot})\right) 
, \\ 
& \Delta_x=
-\tilde{\varepsilon}^2\partial_{z}^2\ttf\,\partial_{\tau}
+\tilde{\varepsilon}^2|\partial_{z}\ttf|^2 \partial_{\tau}^2 
- 2 \tilde{\varepsilon}^2\partial_{z}\ttf \partial_{z}\partial_{\tau}  
+ \tilde{\varepsilon}^2 \partial_z^2 
-\tilde{\varepsilon}\Delta_{\zbot}\ttf \, \partial_{\tau}
+\tilde{\varepsilon}|\nabla_{\zbot}\ttf|^2 \partial_{\tau}^2
-2 \tilde{\varepsilon} \nabla_{\zbot}\ttf\cdot\nabla_{\zbot}\partial_{\tau}
+ \tilde{\varepsilon}\Delta_{\zbot}\\
&\hspace*{0.5cm}=
\frac{1}{c^2} \partial_{\tau}^2 + \tilde{\varepsilon}\Bigl(
-(\tilde{\varepsilon}\partial_z^2\ttf+\Delta_{\zbot}\ttf)\, \partial_{\tau}
-2(\tilde{\varepsilon}\partial_z\ttf\, \partial_z +\nabla_{\zbot}\ttf\cdot\nabla_{\zbot})\, \partial_{\tau}
+\tilde{\varepsilon} \partial_z^2 + \Delta_{\zbot}\Bigr)
\end{aligned}
\end{equation}
Here and in \eqref{parax} we have used the fact that $\partial_t\check{\slo}=0$, $\partial_t\check{\ttf}=0$, hence $\partial_{\tau}\slo=0$, $\partial_{\tau}\ttf=0$, 
$\partial_{x_1}\check{\slo}=\tilde{\varepsilon}\partial_z\slo$, 
$\partial_{x_1}\check{\ttf}=\tilde{\varepsilon}\partial_z\ttf$, 
$\nabla_{x'}\check{\ttf}=\sqrt{\tilde{\varepsilon}}\nabla_{\zbot}\ttf$
}

\medskip



\subsection*{Frequency domain formulation}
We make the time harmonic ansatz 
\begin{align*}
	\phi_k(\tau, z, y) &= \hat{\phi}_k(z,y) e^{\imath\omega_k \tau}  \qquad \text{ for } k \in \left\lbrace 1,2\right\rbrace, \\
	\psi(\tau, z,y)&=\Re{\left\lbrace \hat{\psi}(z,y) e^{\imath(\omega_1-\omega_2)\tau}\right\rbrace}.
\end{align*}
Here the excitation frequencies are assumed to satisfy $\omega_1>\omega_2>0$, with a comparably small difference frequency $\omegad:=\omega_1-\omega_2>0$, and we abbreviate $\omegap:=\omega_1\omega_2\omegad$.
(Typical values in ultrasonics are $\omega_k\sim 2 \,$MHz, $\omegad\sim 50\,$KHz.) 

\medskip

{
Altogether, for $\tilde{\psi}(t,x_1,x')=\psi(\tau,z,\zbot)
=\Re\left(\psih(z,\zbot)e^{\imath\omega_d \tau}\right)$ and 
$\tilde{\psi}(t,x_1,x')
=\Re\left({\check{\psi}}(x_1,x') e^{\imath\omega_d t}\right)$
with $e^{\imath\omega_d t}= e^{\imath\omega_d x_1/c}e^{\imath\omega_d \tau}$
we end up with the transformation in particular of the wave operator
\begin{equation}\label{parax_deriv_freq}
\begin{aligned}
&\Re\left(\left[-\tfrac{1}{c^2}\omega_d^2 -\Delta_x\right]{\check{\psi}}(x_1,x')\,
e^{\imath\omega_d \ttf}\,e^{\imath\omega_d \tau} \right)
=[\tfrac{1}{c^2}\partial_t^2 -\Delta_x]\tilde{\psi}(t,x_1,x')\\
&= \tilde{\varepsilon}\Re\left(
\Bigl[\imath\omega_d (\tilde{\varepsilon}\partial_z^2\ttf+\Delta_{\zbot}\ttf)
+2\imath\omega_d (\tilde{\varepsilon}\partial_z\ttf\, \partial_z +\nabla_{\zbot}\ttf\cdot\nabla_{\zbot})
-\tilde{\varepsilon} \partial_z^2 - \Delta_{\zbot} \Bigr]\psih(z,\zbot)e^{\imath\omega_d \tau}\right)
\end{aligned}
\end{equation}
and likewise for $\tilde{\phi}_k$, $\phih_k$.\\
Due to the boundary condition $\check{\ttf}(0,x')=0, \ x'\in \tilde{\Omega}_{\zbot}$, we also have 
$\nabla_{\zbot}\ttf(0,x')$ and $\Delta_{\zbot}\ttf(0,x')$, which allows us to neglect the terms containing
$\nabla_{\zbot}\ttf$ and $\Delta_{\zbot}\ttf$ to arrive at
\begin{equation}\label{parax_deriv_freq_approx}
\begin{aligned}
\left[-\tfrac{1}{c^2}\omega_d^2 -\Delta_x\right]{\check{\psi}}(x_1,x')\,
e^{\imath\omega_d \ttf}
\approx
\tilde{\varepsilon}\Bigl[\imath\omega_d \partial_z\slo
+2\imath\omega_d \slo\, \partial_z 
-\tilde{\varepsilon} \partial_z^2 - \Delta_{\zbot} \Bigr]\psih(z,\zbot)
\end{aligned}
\end{equation}
}

The outcome of the nonlinear term is more complicated and depends on the interplay between the asymptotic ansatz \eqref{eq:dervib0} governed by $\varepsilon$ and the paraxial approximation \eqref{parax} governed by $\tilde{\varepsilon}$.
Depending on the relation between $\varepsilon$ and $\tilde{\varepsilon}$, this leads to several different models; see \cite{Part_linIP} for more details.

\medskip 

\subsection*{The considered models}
The inverse problem of combined nonlinearity imaging and spound speed reconstruction will be considered for one of these models, namely 
\begin{equation}\label{PDEs_intro}
\begin{aligned}
&{\imath\omega_k \,\slotil\,\phih_k+} 2\imath\omega_k\slo\partial_z\phih_k-\Delta_\zbot \phih_k = 0\mbox{ in }\tilde{\Omega}, \qquad  k\in\{1,2\},\\
&{\imath\omega_d \,\slotil\,\psih+}
2 \imath \omegad \slo \partial_z\psih - \tilde{\varepsilon} \partial_z^2\psih - \Delta_\zbot \psih = 
\imath\omegap \tilde{\varepsilon}^{-1}\,\nlc \,\phih_1\, \overline{\phih}_2\mbox{ in }\tilde{\Omega}
\end{aligned}
\end{equation}
with boundary conditions
\begin{equation}\label{BCs_intro}
\begin{aligned}
&\phih_k(0,\zbot)=h_k(\zbot), \quad \zbot\in\tilde{\Omega}_\zbot\,,\\
&\partial_{\nu_\zbot} \phih_k(z,\zbot)=
-\imath\omega_k \sigma_k \phih_k(z,\zbot), \quad (z,\zbot) \in(0,L)\times\partial\tilde{\Omega}_\zbot,
&&\quad k\in\{1,2\},\\
&
-\partial_z \psih(0,\zbot)=-\imath\omegad \sigma_0\psih(0,\zbot), \quad
\partial_z \psih(L,\zbot)=-\imath\omegad \sigma_L\psih(L,\zbot), \quad \zbot\in \tilde{\Omega}_\zbot,\\
&\partial_{\nu_\zbot} \psih(z,\zbot)=
-\imath\omegad \sigma \psih (z,\zbot), \quad (z,\zbot) \in(0,L)\times\partial\tilde{\Omega}_\zbot,
\end{aligned}
\end{equation}
and $\tilde{\Omega}=(0,L)\times\tilde{\Omega}_\zbot$.

The coefficients we are interested in are related to the sound speed and the nonlinearity parameter by 
\begin{equation}
	\slo=\frac{1}{c}, \qquad \nlc= 
\frac{B/A + 2}{c^4}, 
\end{equation}

In \eqref{BCs_intro}, $h_k$, $k=1,2$ models excitation on part of the boundary and the remaining boundary conditions are supposed to prevent spurious reflections on the boundary of the computational domain $\tilde{\Omega}$.

\begin{remark}\label{rem:excitation}
To model excitation by an array of piezoelectric transducers, continuity of the normal velocity over the transducer-fluid interface would induce Neumann boundary conditions on the velocity potential. In our setting these would read as  
$-\partial_z\phih_k(0,\cdot)=h_k \mbox{ in }\tilde{\Omega}_\zbot$. 
In an analysis analogous to the one we provided here, this could be tackled by differentiating $\phih_k$ with respect to $z$ and considering the corresponding initial value problems for $\phih_k':=\partial_z\phih_k$.
Dirichlet conditions as used here can be justified a surface force -- pressure continuity $F =[\sigma]\cdot \nu = p\nu$ together with the relation $p_k=\rho\partial_t\phi_k$ derived from momentum balance. We refer to, e.g., \cite[eqs (12),(13)]{FlemischKaltenbacherWohlmuth:2006} and the references therein for the proper choice of interface conditions in structure-acoustics coupling.
\end{remark}

Well-posedness of the forward problem \eqref{PDEs_intro}, \eqref{BCs_intro} will be studied in Section~\ref{sec:forward}.

From the point of view of the outgoing wave described by $\hat{\psi}$, the parameter $\tilde{\varepsilon}$ is not necessarily small, since propagation of sound is non-directed for the $\psi$ field. Therefore, alternatively to \eqref{PDEs_intro}, \eqref{BCs_intro}, we will consider 
\begin{equation}\label{PDEsTransf_intro}
\begin{aligned}
&{\imath\omega_k\slotil\,\phih_k+}2\imath\omega_k\slo\partial_z\phih_k-\Delta_\zbot \phih_k = 0\mbox{ in }\tilde{\Omega}, \qquad  k\in\{1,2\},\\
&-\omegad^2 \slo^2 \check{\psi} -\imath\omegad b\Delta_x \check{\psi}- \Delta_x \check{\psi} = \imath\omegap \, 1_{\Omega}\Bigl(\check{\nlc}\,  
P^{-1}\bigl(\phih_1\, \overline{\phih}_2\bigr)\Bigr)\mbox{ in }\Omega
\end{aligned}
\end{equation}
for $\check{\psi}=P^{-1}\psih$, $\check{\nlc}={e^{-\imath\omega_d \ttf}} P^{-1}\nlc$ (where $P$ is defined by $(P\check{\psi})(z,\zbot)=\check{\psi}(x_1,x')$, cf. \eqref{parax_deriv})
 with boundary conditions
\begin{equation}\label{BCsTransf_Intro}
\begin{aligned}
&\phih_k(0,\cdot)=h_k \mbox{ in }\tilde{\Omega}_\zbot\,,
\quad \partial_{\nu_\zbot} \phih_k(z,\zbot)=
-\imath\omega_k \sigma_k \phih_k & \mbox{ in }(0,L)\times\partial\tilde{\Omega}_\zbot
\quad k\in\{1,2\}\\
&\partial_\nu \check{\psi}=-\imath\omegad \sigma \check{\psi} {-\beta \check{\psi}} \mbox{ on }\partial\Omega 
\end{aligned}
\end{equation}
see Section~\ref{sec:forwardHelmholtz} for its well-posedness analysis.

{
Due to \eqref{parax_deriv}, the coefficients in the boundary conditions \eqref{BCs_intro}, \eqref{BCsTransf_Intro} are related by 
\[
\sigma_0(\zbot)=\sigma(0,\zbot)-\slo(0,\zbot), \quad 
\sigma_L(\zbot)=\sigma(L,\zbot)-\slo(L,\zbot).
\]
Thus with the typical choice $\sigma=\slo$, in the propagation direction $z||x_1$, the absorbing boundary condition in the original coordinates becomes a Neumann condition in paraxial coordinates.
}

\medskip

We wish to mention that the paraxial approximation is also made use of in the derivation of the Khokhlov-Zabolotskaya-Kuznetsov (KZK) equation
\begin{equation}\label{KZK}
2\slo \partial^2_{\tau z} \psi - \nabla_y^2 \psi =
\frac{\nlc}{2}\partial_{\tau}\left| \partial_{\tau} \psi\right|^2
\end{equation}
\cite{ZabolotskayaKhokhlov69}, see also, e.g., \cite{Rozanova:2007} for its analysis. Note however, that in our case, the quadratic nonlinearity is decoupled and appears as a source term  for the $\psi$ equation, whereas \eqref{KZK} is a nonlinear (more precisely, quasilinear) equation, whose expansion in frequency domain would lead to an inifinite system of space-dependent PDEs, similarly to \cite{periodicWestervelt}.

\subsection*{The inverse problem}
Our aim is to reconstruct $\slo$ and $\nlc$ in the boundary value problem \eqref{PDEs}, \eqref{BCs}
from measurements of the acoustic pressure 
\begin{equation}\label{obs}
\pmeas=\imath \omegad\mbox{tr}_\Gamma\psih \quad \mbox{ in }\Gamma
\end{equation}
where $\Gamma=P(\check{\Gamma})$ and $\check{\Gamma}$ is a manifold representing the receiver array immersed in the acoustic domain $\Omega$. 

This can be formulated as an operator equation
\begin{equation}\label{Fqy}
F(\slo,\nlc)=y
\end{equation}
with the forward operator $F=C\circ S$ being a concatenation of the (nonlinear) parameter-to-state map 
\begin{equation}\label{par2state}
S:(\slo,\nlc)\mapsto (\phih_1,\phih_2,\psih) \mbox{ solving \eqref{PDEs_intro}, \eqref{BCs_intro}}
\end{equation}
with the (linear) observation operator
\begin{equation}\label{obsop}
C:(\phih_1,\phih_2,\psih)\mapsto \imath \omegad\mbox{tr}_\Gamma\psih .
\end{equation}
We consider $F$ as an operator $F:\Xslo\times\Xnlc\to Y$ with the parameter space $X:=\Xslo\times\Xnlc$ and the data space $Y$ yet to be specified.

In Section~\ref{inverse}, we will alternatively consider an all-at-once formulation of this problem, that keeps the parameters and the state as simultaneous unknowns, thus avoiding the use of a parameter-to-state map $S$.

\medskip

Identification of the nonlinearity coefficient $\eta$ in time domain models of nonlinear acoustics has been studied, e.g., in \cite{AcostaUhlmannZhai2022, nonlinearityimaging_Westervelt, nonlinearityimaging_fracWest, nonlinearityimaging_both, nonlinearityimaging}; in particular \cite{nonlinearityimaging_both} also considers simultaneous identification of $\slo$ and $\nlc$. However, the physical background and therefore also the model differs from the ones we consider here.

A preliminary analysis of the inverse problem of ultrasound vibro-acoustography with models similar to those considered here, but still without the paraxial approximation,  can be found in \cite{vibroacoustics}. In \cite{Part_linIP}, models using a paraxial approximation are derived in time and frequency domain, and the inverse problem of reconstructing $\nlc$ is studied.

\medskip 

The plan of this paper is as follows. Section~\ref{sec:forward} is devoted to the forward problem of solving the PDE model for given coefficients. We  prove well-definedness and differentiability of the parameter-to-state map \eqref{par2state} for both options \eqref{PDEs_intro} and \eqref{PDEsTransf_intro}, in order to justify the application of Newton type methods for the inverse problems. These are discussed in Section~\ref{inverse}, where we establish convergence of a frozen Newton method for reconstructing $\nlc$ and $\slo$ in \eqref{PDEsTransf_intro}. We point to the fact that proving convergence of iterative regularization methods is notoriously difficult in inverse problems with boundary observations, as typical for tomographic imaging. This is due to the fact that convergence criteria, such as the so-called tangential cone condition, can usually not be verified in the situation of restricted observations. We therefore here work with a range invariance condition instead, that indeed can be rigorously verified here and allows to conclude convergence.

\medskip 

An implementation and numerical experiments with the methods analyzed here is subject of future work.
Some numerical results on the simultaneous reconstruction of $\slo$ and $\nlc$ in the Westervelt equation (thus related to this work, but using a model in time domain with a single PDE rather than a system) based on a Newton type iteration as well, can be found in \cite{nonlinearityimaging_both}.

\section{Well-posedness of the forward problem}\label{sec:forward}
In this section we will prove well-definedness and differentiability of the parameter-to-state map for the systems \eqref{PDEs_intro}, \eqref{BCs_intro} and \eqref{PDEsTransf_intro}, \eqref{BCsTransf_Intro}, respectively.

In doing so, we put a particular emphasis on monitoring the smoothness assumptions on the coefficients, which we aim at keeping minimal in view of the fact that in practice, $\slo$ and $\nlc$ tend to be only piecewise smooth and also the inverse problem becomes more ill-posed the higher order the norm in preimage space needs to be chosen.

To justify the use of Newton's method for solving the inverse problem \eqref{Fqy}, we will also prove differentiability of the parameter-to-state map \eqref{par2state}. 

\subsection{Well-posedness of paraxial wave propagation with variable coefficients}
Consider
\begin{equation}\label{PDEs}
\begin{aligned}
&{\imath\omega_k\slotil\,\phih_k+}
2\imath\omega_k\slo\partial_z\phih_k-\Delta_\zbot \phih_k = 0\mbox{ in }\tilde{\Omega}, \qquad  k\in\{1,2\},\\
&{\imath\omegad\slotil\,\psih+}
2\imath \omegad \slo \partial_z\psih - \tilde{\varepsilon} \partial_z^2\psih - \Delta_\zbot \psih = 
\imath\omegap \tilde{\varepsilon}^{-1}\,\nlc \,\phih_1\, \overline{\phih}_2\mbox{ in }\tilde{\Omega}
\end{aligned}
\end{equation}
with boundary conditions
\begin{equation}\label{BCs}
\begin{aligned}
&\phih_k(0,\zbot)=h_k(\zbot), \quad \zbot\in\tilde{\Omega}_\zbot\,,\\
&\partial_{\nu_\zbot} \phih_k(z,\zbot)=
-\imath\omega_k \sigma_k \phih_k(z,\zbot), \quad (z,\zbot) \in(0,L)\times\partial\tilde{\Omega}_\zbot,
&&\quad k\in\{1,2\}\\
&
-\partial_z \psih(0,\zbot)=-\imath\omegad \sigma_0\psih(0,\zbot), \quad
\partial_z \psih(L,\zbot)=-\imath\omegad \sigma_L\psih(L,\zbot), \quad \zbot\in \tilde{\Omega}_\zbot,\\
&\partial_{\nu_\zbot} \psih(z,\zbot)=
-\imath \omegad \sigma \psih (z,\zbot), \quad (z,\zbot) \in(0,L)\times\partial\tilde{\Omega}_\zbot
\end{aligned}
\end{equation}
and $\tilde{\Omega}=(0,L)\times\tilde{\Omega}_\zbot$.

The equations take the form of a (perturbed) Schr\"odinger equation, hence well-known techniques for that equation can be adopted here (see, e.g., \cite{Kenig_lecturenotes:2007} and the references therein).
Since we here require estimates that are explicit in terms of appropriate norms of $\slo$ and $\nlc$, we will first of all provide some energy estimates for solutions of the linear variable coefficient problems  
\begin{equation}\label{Schroedinger}
\begin{aligned}
&\imath\omega \,\slowtil\,u+
2\imath\omega\slo\partial_z u -\Delta_\zbot u = f
\quad \textup{ in }(0,L)\times\tilde{\Omega}_\zbot,\\
&u(0,\zbot)=h(\zbot), \quad \zbot\in\tilde{\Omega}_\zbot\,,\\
&\partial_{\nu_\zbot} u(z,\zbot)=
-\imath\omega\sigma u(z,\zbot), \quad (z,\zbot) \in(0,L)\times\partial\tilde{\Omega}_\zbot
\end{aligned}
\end{equation}
and 
\begin{equation}\label{perturbedSchroedinger}
\begin{aligned}
&{\imath\omega \,\slowtil\,u^b+}
2\imath\omega\slo\partial_z u^b - b\partial_z^2 u^b -\Delta_\zbot u^b = f
\quad \textup{ in }(0,L)\times\tilde{\Omega}_\zbot,\\
&-\partial_z u^b(0,\zbot)=-\imath\omega \sigma_0 u^b(0,\zbot), \quad
\partial_z u^b(L,\zbot)=-\imath\omega \sigma_L u^b(L,\zbot), \quad \zbot\in \tilde{\Omega}_\zbot,\\
&\partial_{\nu_\zbot} u^b(z,\zbot)=
-\imath\omega\sigma u^b(z,\zbot), \quad (z,\zbot) \in(0,L)\times\partial\tilde{\Omega}_\zbot.
\end{aligned}
\end{equation}
Here we assume all coefficients $\slo$, $\sigma_0$, $\sigma_L$, $\sigma$, $b$ to be real valued and positive with 
$\slo,\,\frac{1}{\slo}\,\partial_z\slo\,\in L^\infty(0,L;L^\infty(\tilde{\Omega}_\zbot))$, 
$\sigma,\,\frac{1}{\sigma}\in L^\infty(0,L;L^\infty(\partial\tilde{\Omega}_\zbot))$, 
$\sigma_0,\,\frac{1}{\sigma_0},\,\sigma_L,\,\frac{1}{\sigma_L}\,\in L^\infty(\tilde{\Omega}_\zbot)$
and $b$ constant.
The following two identities will be used repeatedly
\[
2\slo\Re(\partial_zv\, \overline{v}) =
\slo\frac{d}{dz}|v|^2= \frac{d}{dz}\Bigl(\slo |v|^2\Bigr)-\partial_z\slo\,|v|^2,
\]
\[
\int_{\tilde{\Omega}_\zbot} -\Delta_\zbot u\, w\, d\zbot
=\int_{\tilde{\Omega}_\zbot} \nabla_\zbot u\cdot\nabla_\zbot w\, d\zbot
+\imath\omega\int_{\partial\tilde{\Omega}_\zbot}\sigma u\, w \, d\Gamma(\zbot)
\,.
\]
Moreover, for \eqref{perturbedSchroedinger}, we will assume the Poincar\'{e}-Friedrichs type estimate on the domain $\tilde{\Omega}_\zbot$ 
\begin{equation}\label{PF_perturbedSchr}
\|v\|_{\Ltwotil}^2\leq C_0\|\sqrt{\sigma}v\|_{\Ltwotilbdy}^2+C_1\|\nabla_\zbot v\|_{\Ltwotil}^2 
\qquad \forall v\in H^1(\tilde{\Omega}_\zbot)
\end{equation}
to hold for some $C_0$, $C_1>0$.

\begin{proposition}\label{prop:wellposed}
If $\slo \in W^{1,\infty}(0,L;L^\infty(\tilde{\Omega}_\zbot))$, 
{$\partial_z\slowtil\in L^\infty(\tilde{\Omega})$}, 
$h,\, \Delta_\zbot h+f(0)\in\Ltwotil$, $f\in L^\infty(0,L;\Ltwotil)\cap W^{1,1}(0,L;\Ltwotil)$ 
then a solution of \eqref{Schroedinger}, exists, is unique and satisfies the estimate
\[
\begin{aligned}
&\|\nabla_\zbot u(z)\|_{L^\infty(0,L;\Ltwotil)}^2
+\omega\|\sqrt{\slo} \partial_zu\|_{L^\infty(0,L;\Ltwotil)}^2
+\omega\|\sqrt{\slo} u\|_{L^\infty(0,L;\Ltwotil)}^2\\
&\leq
C e^{L(\slo_\sim+\mu)}\Bigl(\omega|\sqrt{\slo(0)} h|_{\Ltwotil}^2
+|\frac{1}{2\sqrt{\slo(0)}}(\Delta_\zbot h+f(0))|_{\Ltwotil}^2
\\&\hspace*{3cm}
+\tfrac{1}{\omega}\left\| \tfrac{1}{\sqrt{\slo}} f\right\|_{L^\infty(0,L;\Ltwotil)}^2
+\left\|\tfrac{1}{\sqrt{\slo}} f\right\|_{L^1(0,L;\Ltwotil)}^2+\left\|\tfrac{1}{\sqrt{\slo}} \partial_z f\right\|_{L^1(0,L;\Ltwotil)}^2\Bigr)
\end{aligned}
\]
with some constants $\mu>0$, $C>0$ independent of $\omega$.

If $\slo \in W^{1,\infty}(0,L;L^\infty(\tilde{\Omega}_\zbot))$, $f\in L^2(0,L;\Ltwotil)$, 
\eqref{PF_perturbedSchr} holds
and the relative smallness conditions on $\slo$, $\partial_z\slo$, $\omega$
\begin{equation}\label{smallness_perturbedSchr}
C_0|{\slowtil-}\partial_z\slo|_{L^\infty} <1 \textup{ and }
2C_1\omega^2|\slo|_{L^\infty}^2< b(1-C_0|{\slowtil-}\partial_z\slo|_{L^\infty})
\end{equation}
are satisfied with $\underline{c}$, $\lambda$ as in \eqref{PF_perturbedSchr},
then a solution of \eqref{perturbedSchroedinger}, exists, is unique and satisfies the estimate
\[
\|u^b\|_{H^1((0,L)\times\tilde{\Omega}_\zbot)}^2\leq C\left(1+\frac{1}{\omega^2}\right) \|f\|_{L^2(0,L;L^2(\Omega)}^2
\]
with some constant $C>0$ independent of $\omega$.
\end{proposition}

\begin{proof}
As in standard evolutionary PDEs, the proof is based on a Galerkin approximation by eigenfunctions of the negative Laplacian, energy estimates and taking weak limits, cf. e.g., \cite{Evans:1998}.
We will here focus on the energy estimates, since the other steps of the proof are relatively straightforward for the linear problems under consideration.

We will multiply the PDEs with appropriate test functions and integrate over $\tilde{\Omega}_\zbot$, abbreviating by $|\cdot|_{\Ltwotil}$, $|\cdot|_{\Ltwotilbdy}$ the $L^2$ norm on $\tilde{\Omega}_\zbot$ and $\partial\tilde{\Omega}_\zbot$, respectively.

Testing this way first of all \eqref{Schroedinger} with $\overline{u}$ and taking the real and imaginary parts we obtain
\begin{equation}\label{enid1Re}
|\nabla_\zbot u(z)|_{\Ltwotil}^2=2\omega\int_{\tilde{\Omega}_\zbot} \slo\Im(\partial_zu(z)\, \overline{u}(z)))\, d\zbot 
+\int_{\tilde{\Omega}_\zbot} \Re(f\, \overline{u})\, d\zbot,
\end{equation}
\begin{equation}\label{enid1Im}
\begin{aligned}
\omega\frac{d}{dz} |\sqrt{\slo} u|_{\Ltwotil}^2(z)
+\omega|\sqrt{\sigma}u(z)|_{\Ltwotilbdy}^2 &=-\omega\int_{\tilde{\Omega}_\zbot}({\slowtil-}\partial_z\slo)\,|u(z)|^2\, d\zbot
+\int_{\tilde{\Omega}_\zbot} \Im(f\, \overline{u})\, d\zbot\\
&\leq \omega(\slo_\sim+\mu) |\sqrt{\slo} u|_{\Ltwotil}^2(z) + \frac{1}{4\mu\omega}|\frac{1}{\sqrt{\slo}} f|_{\Ltwotil}^2(z)
\end{aligned}
\end{equation}
for all $z\in(0,L)$, with $\slo_\sim:=\|\frac{\slowtil-\partial_z\slo}{\slo}\|_{L^\infty(0,L;L^\infty(\tilde{\Omega}_y))}$, where we have used Young's inequality.


Doing the same with \eqref{perturbedSchroedinger} and integrating over $(0,L)$ we obtain
\begin{equation}\label{enid1Re_perturbed}
\begin{aligned}
&\|\nabla_\zbot u^b\|_{L^2(0,L;\Ltwotil)}^2
+b\|\partial_z u^b\|_{L^2(0,L;\Ltwotil)}^2
=2\omega\int_0^L\int_{\tilde{\Omega}_\zbot} \slo\Im(\partial_zu^b\, \overline{u}^b)\, d\zbot \, dz
+\int_0^L\int_{\tilde{\Omega}_\zbot} \Re(f\, \overline{u}^b)\, d\zbot\, dz,
\end{aligned}
\end{equation}
\begin{equation}\label{enid1Im_perturbed}
\begin{aligned}
&\omega|\sqrt{\slo} u^b|_{\Ltwotil}^2(L)
+\omega\|\sqrt{\sigma}u^b\|_{L^2(0,L;\Ltwotilbdy)}^2 + b\omega\sigma_0|u^b(0)|_{\Ltwotil}^2+ b\omega\sigma_L|u^b(L)|_{\Ltwotil}^2\\
&=\omega|\sqrt{\slo} u^b|_{\Ltwotil}^2(0)
-\omega\int_0^L\int_{\tilde{\Omega}_\zbot}({\slowtil-}\partial_z\slo)\,|u^b|^2\, d\zbot \, dz
+\int_0^L\int_{\tilde{\Omega}_\zbot} \Im(f\, \overline{u}^b)\, d\zbot \, dz.
\end{aligned}
\end{equation}


We also differentiate \eqref{Schroedinger} with respect to $z$ and test with $\partial_z\overline{u}$. Analogously to above, this yields
\begin{equation}\label{enid2Re}
|\nabla_\zbot \partial_z u(z)|_{\Ltwotil}^2=\omega\int_{\tilde{\Omega}_\zbot} \Bigl(
{\partial_z\slowtil\Im(\partial_zu(z)\, \overline{u}(z)))+}
2\slo\Im(\partial_z^2u(z)\, \partial_z\overline{u}(z)))\Bigr)\, d\zbot 
+\int_{\tilde{\Omega}_\zbot} \Re(\partial_z f\, \partial_z\overline{u})\, d\zbot,
\end{equation}
\begin{equation}\label{enid2Im}
\begin{aligned}
&
{\omega\int_{\tilde{\Omega}_y}\partial_z\slowtil\Re(\partial_zu(z)\, \overline{u}(z)))\,d\zbot+}
\omega\frac{d}{dz} |\sqrt{\slo} \partial_zu|_{\Ltwotil}^2(z)
+\omega|\sqrt{\sigma}\partial_zu(z)|_{\Ltwotilbdy}^2\\ 
&\leq \omega(\slo_+ +\mu) |\sqrt{\slo} \partial_zu|_{\Ltwotil}^2(z) + \frac{1}{4\mu\omega}|\frac{1}{\sqrt{\slo}} \partial_zf|_{\Ltwotil}^2(z).
\end{aligned}
\end{equation}
with $\slo_+= \|\frac{\slowtil+\partial_z\slo}{\slo}\|_{L^\infty(0,L;L^\infty(\tilde{\Omega}_y))}$.

For the initial value problem (with respect to $z$) for the Schr\"odinger equation \eqref{Schroedinger}, we combine {\eqref{enid1Im}+$\rho\cdot$\eqref{enid2Im} yields, for 
$\eta(z)=\eta_0(z)+\rho\eta_1(z)=|\sqrt{\slo} u|_{\Ltwotil}^2(z)+
\rho|\sqrt{\slo} \partial_zu|_{\Ltwotil}^2(z)$ 
that
\[
\begin{aligned}
\frac{d}{dz}\eta_0 &\leq (\slo_\sim+\mu) \eta_0 +\frac{1}{4\omega^2\mu} |\frac{1}{\sqrt{\slo}} f|_{\Ltwotil}^2\\
\frac{d}{dz}\eta_1 &\leq (\slo_+ +\mu) \eta_1 +\frac{1}{4\omega^2\mu} |\frac{1}{\sqrt{\slo}} \partial_zf|_{\Ltwotil}^2
+\|\partial_z\slowtil/\slo\|_{L^\infty(\tilde{\Omega})}\sqrt{\eta_0\eta_1}
\end{aligned}
\]
thus, using Young's inequality with factors $\frac{1}{2\sqrt{\rho}}$, $\frac{\sqrt{\rho}}{2}$ for the last term 
\[
\begin{aligned}
\frac{d}{dz}\eta \leq \beta \eta + \alpha \textup{ where } &
\beta=\max\{\slo_\sim,\slo_+\}+\mu+\sqrt{\rho}/2\,\|\partial_z\slowtil/\slo\|_{L^\infty(\tilde{\Omega})},\\
&\alpha(z)=\frac{1}{4\omega^2\mu}\Bigl(|\frac{1}{\sqrt{\slo}} f(z)|_{\Ltwotil}^2+|\frac{1}{\sqrt{\slo}} \partial_zf(z)|_{\Ltwotil}^2\Bigr).
\end{aligned}
\]
An application of Gronwall's inequality yields
\[
\eta(z)=|\sqrt{\slo} u|_{\Ltwotil}^2(z)+\rho |\sqrt{\slo} \partial_zu|_{\Ltwotil}^2(z)
\leq e^{\beta z}
\Bigl(|\sqrt{\slo} u|_{\Ltwotil}^2(0)+\rho|\sqrt{\slo} \partial_zu|_{\Ltwotil}^2(0)\Bigr)
\int_0^z e^{\beta(z-a)}\alpha(a)\, da,
\]
}
where we can insert the initial conditions and their differentiated version substituting from the PDE
\[
\sqrt{\slo(0)}u(0)=\sqrt{\slo(0)}h, \qquad
\sqrt{\slo(0)}\partial_z u(0)  = \frac{1}{2\imath\omega\sqrt{\slo(0)}}\Bigl(f(0) +\Delta_\zbot h -\imath\omega \,\slowtil(0)\,h\Bigr).
\]

Using this together with the Cauchy-Schwarz and Young's inequalities in \eqref{enid1Re} yields the estimate
\begin{equation}\label{enest0}
\begin{aligned}
&|\nabla_\zbot u(z)|_{\Ltwotil}^2
\leq 2\omega\int_{\tilde{\Omega}_\zbot} |\sqrt{\slo} \partial_zu|(z)\,
|\sqrt{\slo} u|(z)\, d\zbot
+\int_{\tilde{\Omega}_\zbot} |\tfrac{1}{\sqrt{\slo}} f|\, |\sqrt{\slo} u|\, d\zbot\\
&\leq \omega\lambda_1|\sqrt{\slo} \partial_zu|_{\Ltwotil}^2(z) 
+\tfrac{\omega}{\lambda_1}|\sqrt{\slo} u|_{\Ltwotil}^2(z)
+\tfrac{\lambda_2}{2}|\sqrt{\slo} u|_{\Ltwotil}^2(z)
+\tfrac{1}{2\lambda_2}|\tfrac{1}{\sqrt{\slo}} f|_{\Ltwotil}^2(z)\\
&{\leq \omega\lambda_1\eta + \tfrac{\rho\lambda_1}{4\omega(\lambda_1^2-\rho)}|\tfrac{1}{\sqrt{\slo}} f|_{\Ltwotil}^2(z)}
\end{aligned}
\end{equation}
{
where we have chosen $\lambda_2$ such that $\omega\lambda_1=\rho(\tfrac{\omega}{\lambda_1}+\tfrac{\lambda_2}{2})$ in order to make optimal use of $\eta$.
}

\medskip

To obtain an estimate for the two point boundary value problem (with respect to $z$) for the perturbed Schr\"odinger equation \eqref{perturbedSchroedinger} we combine \eqref{enid1Re_perturbed} $+\frac{C_0}{C_1\omega}\cdot$ \eqref{enid1Im_perturbed} and use \eqref{PF_perturbedSchr}, \eqref{smallness_perturbedSchr} to obtain
\begin{equation*}
\begin{aligned}
&\frac{1}{C_1} \|u^b\|_{L^2(0,L;\Ltwotil)}^2+b\|\partial_z u^b\|_{L^2(0,L;\Ltwotil)}^2\\
&\leq\|\nabla_\zbot u^b\|_{L^2(0,L;\Ltwotil)}^2
+b\|\partial_z u^b\|_{L^2(0,L;\Ltwotil)}^2
+\frac{C_0}{C_1}|\sqrt{\slo} u^b|_{\Ltwotil}^2(L)\\
&+\frac{C_0}{C_1}\|\sqrt{\sigma}u^b\|_{L^2(0,L;\Ltwotilbdy)}^2
+ \frac{C_0}{C_1}(b\sigma_0-|\slo|_{L^\infty})|u^b(0)|_{\Ltwotil}^2
+ \frac{C_0}{C_1} b\sigma_L|u^b(L)|_{\Ltwotil}^2\\
&\leq 2\omega\int_0^L\int_{\tilde{\Omega}_\zbot} \slo|\Im(\partial_zu^b\, \overline{u}^b)|\, d\zbot \, dz
+\int_0^L\int_{\tilde{\Omega}_\zbot} |\Re(f\, \overline{u}^b)|\, d\zbot\, dz
-\frac{C_0}{C_1} \int_0^L\int_{\tilde{\Omega}_\zbot} ({\slowtil-}\partial_z\slo)\,|u^b|^2\, d\zbot \, dz
+\frac{C_0}{C_1\omega} \int_0^L\int_{\tilde{\Omega}_\zbot} |\Im(f\, \overline{u}^b)|\, d\zbot \, dz
\\
&\leq (\mu+\tfrac{C_0}{C_1}|{\slowtil-}\partial_z\slo|_{L^\infty}) \|u^b\|_{L^2(0,L;\Ltwotil)}^2
\\&\qquad
+\frac{1}{2\mu}\Bigl(
4\omega^2 |\slo|_{L^\infty}^2 
\|\partial_z u^b\|_{L^2(0,L;\Ltwotil)}^2
+(1+(\tfrac{C_0}{C_1\omega})^2)\|f\|_{L^2(0,L;\Ltwotil)}^2 \Bigr),
\end{aligned}
\end{equation*}
where $C_0$, $C_1$ are as in \eqref{PF_perturbedSchr}.
In order to be able to achieve 
\[
C_1\mu+C_0|{\slowtil-}\partial_z\slo|_{L^\infty} <1\quad \textup{ and }\quad
\frac{2\omega^2|\slo|_{L^\infty}^2}{\mu} < b,
\]
we invoke \eqref{smallness_perturbedSchr} and choose 
$\mu\in(
\frac{2\omega^2|\slo|_{L^\infty}^2}{b},
\frac{1-C_0|{\slowtil-}\partial_z\slo|_{L^\infty}}{C_1})$.
\end{proof}

\begin{remark}\label{rem:smallness_omega}
Since we will apply the part of Proposition~\ref{prop:wellposed} concerning \eqref{Schroedinger} with large $\omega=\omega_k$, $k=1,2$, we track frequency dependence of the estimates explicitly there.
On the other hand, the energy estimate in Proposition~\ref{prop:wellposed} for \eqref{perturbedSchroedinger} will only be used for the (small) difference frequency $\omega=\omegad$ and therefore the restrictions on the size of $\omega$ made there do not compromise practical applicability in our context.

{
A wavenumber explicit estimate for \eqref{perturbedSchroedinger} follows by application of the transformation \eqref{parax_deriv_freq} to the Helmholtz equation without having to impose a smallness condition \eqref{smallness_perturbedSchr}, see Proposition~\ref{prop:wellposedHelmholtz} in the next subsection.
}
\end{remark}

Using Proposition~\ref{prop:wellposed}, we can conclude the following results on the parameter-to-state map 
\begin{equation}\label{par2state_secforward}
S:(\slo,\nlc)\mapsto (\phih_1,\phih_2,\psih) \mbox{ solving \eqref{PDEs}, \eqref{BCs}}
\end{equation}
and its linearization $dS(\slo,\nlc)(\dslo,\dnlc):=(\dphi_1,\dphi_2,\dpsi)$ defined by the system
\begin{equation}\label{dPDEs}
\begin{aligned}
&\imath\omega_k\slotil\,\dphi_k+
2\imath\omega_k\slo\partial_z\dphi_k-\Delta_\zbot \dphi_k = 
{-\imath\omega_k\partial_z\dslo\,\phih_k}
-2\imath\omega_k\dslo\partial_z\phih_k
\mbox{ in }\tilde{\Omega}, \qquad  k\in\{1,2\},\\
&\imath\omegad\slotil\,\dpsi+
2 \imath \omegad \slo \partial_z\dpsi - \tilde{\varepsilon} \partial_z^2\dpsi - \Delta_\zbot \dpsi 
= {-\imath\omegad\partial_z\dslo\,\psih}
-2 \imath \omegad \dslo \partial_z\psih\\
&\hspace*{3.5cm}
+\imath\omegap \tilde{\varepsilon}^{-1}\,\Bigl(
\dnlc \,\phih_1\, \overline{\phih}_2
+\nlc (\dphi_1\, \overline{\phih}_2+\phih_1\, \overline{\dphi}_2)
\Bigr) \mbox{ in }\tilde{\Omega}
\end{aligned}
\end{equation}
with boundary conditions
\begin{equation}\label{dBCs}
\begin{aligned}
&\dphi_k(0,\zbot)=0, \quad \zbot\in\tilde{\Omega}_\zbot\,,\\
&\partial_{\nu_\zbot} \dphi_k(z,\zbot)=
-\imath\omega_k \sigma_k \dphi_k(z,\zbot), \quad (z,\zbot) \in(0,L)\times\partial\tilde{\Omega}_\zbot,
&&\quad k\in\{1,2\},\\
&
-\partial_z \dpsi(0,\zbot)=-\imath\omegad\sigma_0\dpsi(0,\zbot), \quad
\partial_z \dpsi(L,\zbot)=-\imath\omegad\sigma_L\dpsi(L,\zbot), \quad \zbot\in \tilde{\Omega}_\zbot,\\
&\partial_{\nu_\zbot} \dpsi(z,\zbot)=
-\imath\omegad\sigma \dpsi (z,\zbot), \quad (z,\zbot) \in(0,L)\times\partial\tilde{\Omega}_\zbot.
\end{aligned}
\end{equation}

\begin{corollary}\label{cor:par2state}
Assume $h,\, \Delta_\zbot h\in\Ltwotil$. 
Then $S$ is well-defined by \eqref{par2state_secforward} as an operator 
\begin{equation*}
S:\mathcal{D}_\slo^1\times L^2(0,L;L^p(\tilde{\Omega}_\zbot))\to V
\end{equation*}
and G\^{a}teaux differentiable for all $(\slo,\nlc)\in\mathcal{D}_\slo^2\times L^2(0,L;L^p(\tilde{\Omega}_\zbot))$ as an operator \\
$W^{1,\infty}(0,L;L^\infty(\tilde{\Omega}_\zbot))\times L^2(0,L;L^p(\tilde{\Omega}_\zbot))\to V$, where
\begin{equation}\label{defD1D2V}
\begin{aligned}
\mathcal{D}_\slo^1&=\{\slo \in W^{1,\infty}(0,L;L^\infty(\tilde{\Omega}_\zbot))\, : \, \textup{ \eqref{smallness_perturbedSchr} holds for } 
\omega=\omegad=\omega_1-\omega_2
\}, \\ 
\mathcal{D}_\slo^2&=\mathcal{D}_\slo^1\cap W^{2,q}(0,L;L^\infty(\tilde{\Omega}_\zbot)),\\
V&=\Bigl(L^\infty(0,L;H^1(\tilde{\Omega}_\zbot)\cap W^{1,\infty}(0,L;L^2(\tilde{\Omega}_\zbot)\Bigr)^2\times H^1((0,L)\times\tilde{\Omega}_\zbot)
\end{aligned}
\end{equation}
for some $p,q\in\begin{cases}
&(1,\infty]\textup{ if $\tilde{\Omega}_\zbot\subseteq\mathbb{R}^2$}\\
&[1,\infty]\textup{ if $\tilde{\Omega}_\zbot\subseteq\mathbb{R}^1$}\end{cases}$, with $dS(\slo,\nlc)(\dslo,\dnlc):=(\dphi_1,\dphi_2,\dpsi)$ defined by \eqref{dPDEs}, \eqref{dBCs}.
\end{corollary}
\begin{proof}
The first part of Proposition~\ref{prop:wellposed} with $\omega=\omega_k$ implies that the solutions $\phih_1,\, \phih_2\,\in L^\infty(0,L;H^1(\tilde{\Omega}_y))\cap W^{1,\infty}(0,L;\Ltwotil)$.\\ 
Setting $f:=\imath\omegap \tilde{\varepsilon}^{-1}\,\nlc \,\phih_1\, \overline{\phih}_2\,\in L^2(0,L;\Ltwotil)$ and since $H^1(\tilde{\Omega}_y)$ continuously embeds into $L^{2p/(p-1)}$, we conclude $\psih\in H^1((0,L)\times\tilde{\Omega}_\zbot)$ from the second part of Proposition~\ref{prop:wellposed} with $\omega=\omegad$.

It is readily checked that the first variation of $S$ at $(\slo,\nlc)$ in direction $(\dslo,\dnlc)$ is given by the solution of \eqref{dPDEs}, \eqref{dBCs} and that $dS(\slo,\nlc)$ is a linear operator. To prove boundedness of $dS(\slo,\nlc):W^{1,\infty}(0,L;L^\infty(\tilde{\Omega}_\zbot))\to V$ for $(\slo,\nlc)\in\mathcal{D}_\slo^2\times L^2(0,L;L^p(\tilde{\Omega}_\zbot))$ 
we apply the first part of Proposition~\ref{prop:wellposed} with $\omega=\omega_k$ and $f:=-2\imath\omega_k\dslo\partial_z\phih_k$, $h=0$. To prove that $f\in L^\infty(0,L;\Ltwotil)\cap W^{1,1}(0,L;\Ltwotil)$ we need somewhat more regularity of $\phih_k$. The latter can be derived by considering the PDE for $\phih_k:=\partial_z\phih_k$, \eqref{Schroedinger} with $f=-2\imath\omega_k\partial_z\slo\partial_z\phih_k\in L^\infty(0,L;\Ltwotil)\cap W^{1,1}(0,L;\Ltwotil)$ for $\slo\in\mathcal{D}_\slo^2$ and implies 
$\phih_k\in W^{1,\infty}(0,L;H^1(\tilde{\Omega}_\zbot))\cap W^{1,\infty}(0,L;\Ltwotil)$.\\
Finally, application of the second part of Proposition~\ref{prop:wellposed} with $\omega=\omegad$ and $f:=-2 \imath \omegad \dslo \partial_z\psih
+\imath\omegap \tilde{\varepsilon}^{-1}\,\Bigl(
\dnlc \,\phih_1\, \overline{\phih}_2
+\nlc (\dphi_1\, \overline{\phih}_2+\phih_1\, \overline{\dphi}_2)
\Bigr)\in L^2(0,L;\Ltwotil)$ with $f(0)\in\Ltwotil$ implies the required bound on $\dpsi$ in $H^1((0,L)\times\tilde{\Omega}_\zbot)$.
\end{proof}

\subsection{Well-posedness, using the Helmholtz equation for the propagating wave}\label{sec:forwardHelmholtz}

Alternatively to the paraxial form of the PDE for $\psih$ we consider
\begin{equation}\label{PDEsTransf}
\begin{aligned}
&{\imath\omega_k\slotil\,\phih_k+}
2\imath\omega_k\slo\partial_z\phih_k-\Delta_\zbot \phih_k = 0\mbox{ in }\tilde{\Omega}, \qquad  k\in\{1,2\},\\
&-\omegad^2 \slo^2 \check{\psi} - \Delta_x \check{\psi} = \imath\omegap \, 1_{\Omega}\Bigl(\check{\nlc}\,  
P^{-1}\bigl(\phih_1\, \overline{\phih}_2\bigr)\Bigr)\mbox{ in }\Omega
\end{aligned}
\end{equation}
for $\check{\psi}=P^{-1}\psih$, $\check{\nlc}={e^{-\imath\omega_d \ttf}} P^{-1}\nlc$ (cf \eqref{parax_deriv}) with boundary conditions
\begin{equation}\label{BCsTransf}
\begin{aligned}
&\phih_k(0,\cdot)=h_k \mbox{ in }\tilde{\Omega}_\zbot\,,
\quad \partial_{\nu_\zbot} \phih_k(z,\zbot)=
-\imath\omega_k\sigma_k \phih_k & \mbox{ in }(0,L)\times\partial\tilde{\Omega}_\zbot,
\quad k\in\{1,2\},\\
&\partial_\nu \psic=-\imath\omegad\sigma \psic  \mbox{ on }\partial\Omega 
\end{aligned}
\end{equation}
with $P(\Omega)\supseteq\tilde{\Omega}=(0,L)\times\tilde{\Omega}_\zbot$, where $P$ is the paraxial transform defined by \eqref{parax} and 
$1_{\Omega}$ the extension by zero operator for a function defined on $P^{-1}(\tilde{\Omega})$ to all of $\Omega$. This allows to use a possibly larger propagation domain for the $\psi$ wave as compared to the one for the $\phi_1$, $\phi_2$ waves, and to get rid of the smallness assumption \eqref{smallness_perturbedSchr}.

Well-posedness and G\^{a}teaux differentiability of the forward operator defined by the system \eqref{PDEsTransf}, \eqref{BCsTransf} follows analogously to Corollary~\ref{cor:par2state} by combining the first part of Proposition~\ref{prop:wellposed} with known results on the Helmholtz equation with impedance boundary conditions, see, e.g., \cite[Chapter 8]{ThesisMelenk}.
For completeness and to track the required coefficient regularity we here provide the essential arguments. 

We start with some energy estimates for the Helmholtz equation that can be obtained by applying the testing strategy from \cite[Chapter 8]{ThesisMelenk} (where the constant coefficient case with $\sigma=\slo$ is considered) to our a variable coefficient setting. Again, these results are in principle available in the literature, see, e.g., \cite{GrahamSauter:2020}, but we aim at making the required regularity of $\slo$ visible for our purposes. 

Testing the general Helmholtz equation on a smooth domain $\Omega$  with impedance boundary conditions 
\begin{equation}\label{genHelmholtz} 
\begin{aligned}
-\omega^2 \slo^2 u - \Delta_x u &= f \mbox{ in }\Omega\\
\partial_\nu u+\imath\omega\sigma u &= h  
\mbox{ on }\partial\Omega 
\end{aligned}
\end{equation}
with $u$ and taking real and imaginary parts yields
\begin{equation}\label{enid_Helmholtz}
\begin{aligned}
&\|\nabla_x u\|_{L^2(\Omega)}^2
-\omega^2\|\slo u\|_{L^2(\Omega)}^2
=\int_\Omega \Re(f\, u)\, dx ,
\\
&\omega\int_{\partial\Omega}\sigma\,|u|^2\, d\Gamma  
=\int_\Omega \Im(f\, u)\, dx.
\end{aligned}
\end{equation}
On a starlike domain $\Omega$ in two space dimensions $d=2$ with 
\begin{equation}\label{cOmega}
x\cdot\nu(x)\geq c_\Omega>0 
\end{equation}
on $\partial\Omega$ we can also use 
$v(x):=x\cdot\nabla u(x)$ as a test function (provided it is contained in $H^1(\Omega)$) and consider the real part to arrive at
\begin{equation}\label{enid_Helmholtz2}
\begin{aligned}
\omega^2\int_\Omega (1+x\cdot\tfrac{\nabla\slo}{\slo})|\slo\, u|^2\, dx 
-\frac12 \omega^2\int_{\partial\Omega} |\slo\, u|^2\, x\cdot\nu(x)\,d\Gamma(x) 
+\frac12 \int_{\partial\Omega} |\nabla u|^2\, x\cdot\nu(x)\,d\Gamma(x) \\
+\omega\Re\left(\imath\sigma \int_{\partial\Omega} u \ x\cdot \overline{\nabla u}\,d\Gamma(x)\right)
=\Re\left(\int_\Omega f\, x\cdot \overline{\nabla u}\,dx\right).
\end{aligned}
\end{equation}
Here we have used the identities 
\[
\begin{aligned}
&2\Re\left(\slo^2 u\, x\cdot\nabla u \right)=
x\cdot \Bigl(\nabla |\slo u|^2 -  2\slo\nabla\slo |u|^2\Bigr),\\ 
&\nabla(x\cdot \overline{\nabla u}) = 
\overline{\nabla u} + (x \cdot\nabla)\overline{\nabla u}, \quad
2\Re\left(\nabla u \cdot (x \cdot\nabla)\overline{\nabla u}\right) = 
\nabla\cdot(|\nabla u|^2 x) - d |\nabla u|^2
\end{aligned}
\]
as well as the divergence Theorem and $d=2$.

This yields the following results in a low and higher regularity regime of $\slo$ and correspondingly estimates on $u$ with different frequency dependence of the constants.
\begin{proposition}\label{prop:wellposedHelmholtz}
If $\Omega$ is a bounded $C^{1,1}$ smooth or convex 
domain, $\slo 
\in L^\infty(\Omega)$, $f\in H^1(\Omega)^*$, 
$\sigma\not=0$ a.e on $\partial\Omega$, then a solution of \eqref{genHelmholtz} exists, is unique and satisfies the estimate
\[
\|\nabla_x u\|_{L^2(\Omega)}^2 + \|u\|_{L^2(\Omega)}^2 \leq 
C(\omega) \bigl( \|f\|_{H^1(\Omega)^*}^2+\|h\|_{H^{-1/2}(\partial\Omega)}^2\bigr)
\]
with some constant $C(\omega)>0$ depending on $\omega$.

If additionally $f\in L^2(\Omega)$, $h\in H^{1/2}(\partial\Omega)$ then
\begin{equation}\label{H2omegadep}
\|D_x^2 u\|_{L^2(\Omega)}^2+\|\nabla_x u\|_{L^2(\Omega)}^2 + \|u\|_{L^2(\Omega)}^2 \leq 
C(\omega) \bigl( \|f\|_{L^2(\Omega)}^2+\|h\|_{H^{1/2}(\partial\Omega)}^2\bigr).
\end{equation}
Here $D_x^2u$ denotes the (weak derivative) Hessian of $u$.

If additionally $\Omega$ is a star-shaped domain with \eqref{cOmega}, $\slo\geq0$, $\ln(\slo) \in W^{1,\infty}(\Omega)$, $f\in L^2(\Omega)$, $h=0$, $d=2$
and the relative smallness/largeness conditions 
\begin{equation}\label{smallness_Helmholtz}
c_\slo=\sup_{x\in\Omega} (-x\cdot\tfrac{\nabla\slo(x)}{\slo(x)})<1, \quad 
\sigma \geq c_\sigma>0
\end{equation}
are satisfied, then the solution of \eqref{genHelmholtz} satisfies the estimate
\begin{equation}\label{estHelmholtz2}
\frac{1}{\omega^2}\|D_x^2 u\|_{L^2(\Omega)}^2+\|\nabla_x u\|_{L^2(\Omega)}^2 + \omega^2 \|u\|_{L^2(\Omega)}^2 
\leq C 
\left(1+\frac{1}{\omega^2}\right) \|f\|_{L^2(\Omega)}^2
\end{equation}
with some $C$ independent of $\omega$. 
\end{proposition}
\begin{proof}
The well-posedness proof in the low regularity regime follows analogously to the one in the constant coefficient case \cite[Chapter 8]{ThesisMelenk}. Since it can hardly be found in the literature for our setting (variably slowness, impedance boundary conditions) we provide it here.
We rewrite \eqref{genHelmholtz} as 
\[
(-\Delta+\text{id}) u = F + Gu
\]
with the bounded and boundedly invertible operator $(-\Delta+\text{id}):H^1(\Omega)\to H^1(\Omega)^*$, the bounded operator $G:H^1(\Omega)\to H^1(\Omega)^*$ and the element $F\in H^1(\Omega)^*$ defined by 
\[
\begin{aligned}
&\langle (-\Delta+\text{id})u,v\rangle_{(H^1)^*,H^1}:= 
\int_\Omega (\nabla u\cdot\nabla\overline{v}+u\,\overline{v})\, dx, \\
&\langle Gu,v\rangle_{(H^1)^*,H^1}:= 
\int_\Omega (1+\omega^2\slo^2)\,u\,\overline{v}\, dx -\imath\omega\int_{\partial\Omega}\sigma u\overline{v}\, d\Gamma,\\ 
&\langle F,v\rangle_{(H^1)^*,H^1}:= \int_\Omega f\,\overline{v}\, dx +\int_{\partial\Omega}h\overline{v}\, d\Gamma
\end{aligned}
\]
for any $v\in H^1(\Omega)$.
This is equivalent to 
\[
(I-K)u = b \quad \text{ with } K = (-\Delta+\text{id})^{-1}G, \ b = (-\Delta+\text{id})^{-1}F,
\]
where $b\in H^1(\Omega)$ and the operator $K$ is bounded from $L^2(\Omega)$ to $H^1(\Omega)$  and thus compact when considered as an operator from $L^2(\Omega)$ into itself.
We can therefore apply Fredholm's alternative for showing bijectivity of $I-K$, then the Bounded Inverse Theorem on the Banach space $L^2(\Omega)$ and finally lift regularity of $u$ to $H^1(\Omega)$ by means of the fixed point identity $u=Ku+b$.
In order to use Fredholm's alternative, we have to prove injectivity of $I-K$. To this end, assume that $(I-K)w=0$ for some $w\in L^2(\Omega)$, which due to $w=Kw$ is automatically contained in $H^1(\Omega)$. The second line of \eqref{enid_Helmholtz} implies that $w$ satisfies homogeneous Dirichlet, hence by the impedance condition also homogeneous Neumann boundary conditions and can therefore be extended by zero to an $H^2(\mathbb{R}^3)$ function\footnote{$w\in H^2(\Omega)$ follows from elliptic regularity and $-\Delta_x w =\omega^2\slo^2 w\in L^2(\Omega)$ provided $\Omega$ is $C^{1,1}$ or convex; This $H^2$ requirement is only needed because of the used uniqueness results from \cite[Section 8.3]{ColtonKress}; 
}
 $\hat{w}$ that then satisfies the homogeneous Helmholtz equation $-\Delta_x \hat{w} +\hat{s}\hat{w}$ on all of $\mathbb{R}^3$ with $\hat{s}:=1_\Omega (s-1) + 1$ and Sommerfeld radiation conditions. From the results in \cite[Section 8.3]{ColtonKress} we conclude that $\hat{w}\equiv0$ and thus $w\equiv0$.

The higher order regularity result \eqref{H2omegadep} follows from elliptic reguarity and the fact that $u$ satisfies $-\Delta_x u = \tilde{f}$, $\partial_\nu u = \tilde{h}$
with 
$\tilde{f}= f+\omega^2 \slo^2 u\in L^2(\Omega)$, 
$\tilde{h}= h-\imath\omega\sigma \textup{tr}_{\partial\Omega}u\in H^{1/2}(\partial\Omega).$  

\medskip

From \eqref{enid_Helmholtz2} with \eqref{cOmega}, \eqref{smallness_Helmholtz}, applying  Cauchy-Schwarz' inequality
we conclude
\begin{equation}
\begin{aligned}
&\omega^2(1-c_\Omega)\|\slo u\|_{L^2(\Omega)}^2 + \frac12\gamma\|\nabla u\|_{L^2(\partial \Omega)}^2 \leq  \frac12 \omega^2\int_{\partial\Omega} |\slo\, u|^2\, x\cdot\nu(x)\,d\Gamma(x)\\ 
&\qquad+\omega^2\Bigl(\sqrt{\int_{\partial\Omega} \sigma^2 |x|^2 |u|^2 \,d\Gamma(x)}\, 
+\sqrt{\int_{\partial\Omega} |x|^2 |f|^2 \,d\Gamma(x)}\Bigr) \|\nabla u\|_{L^2(\partial \Omega)}.
\end{aligned}
\end{equation}
Hence by Young's inequality and the second line in \eqref{enid_Helmholtz}, that is, 
$\omega\|\sqrt{\sigma} u\|_{L^2(\partial\Omega)}^2  \leq \|f/\slo\|_{L^2(\Omega)} \|\slo u\|_{L^2(\Omega)}$,
\[
\begin{aligned}
&\omega^2(1-c_\Omega)\|\slo u\|_{L^2(\Omega)}^2 
\\
&\leq  \frac{\omega^2}{2}\int_{\partial\Omega} |\slo\, u|^2\, x\cdot\nu(x)\,d\Gamma(x) 
+\frac{\omega^2}{\gamma}\int_{\partial\Omega} \sigma^2 |x|^2 |u|^2 \,d\Gamma(x)\, 
+\frac{1}{\gamma}\int_{\partial\Omega} |x|^2 |f|^2 \,d\Gamma(x)\\
&\leq  \omega
\Bigl(\frac{\textup{d}(\Omega)}{2}\,\left\|\frac{\slo^2}{\sigma}\right\|_{L^\infty(\Omega)}
+\frac{\textup{d}(\Omega)^2}{\gamma}\|\sigma\|_{L^\infty(\Omega)}\Bigr) \|f/\slo\|_{L^2(\Omega)} \|\slo u\|_{L^2(\Omega)} \ + \frac{\textup{d}(\Omega)^2}{\gamma} \|f\|_{L^2(\Omega)}^2
\end{aligned}
\]
with $d(\Omega)=\sup_{x\in\Omega}|x|$, hence by another application of Young's inequality
(with factors $\frac{1-c_\Omega}{2}$ and $\frac{1}{2(1-c_\Omega)}$) 
\[
\omega^2\|\slo u\|_{L^2(\Omega)}^2 \leq C \|f\|_{L^2(\Omega)}^2
\]
with 
$C = \frac{1}{(1-c_\Omega)^2} \Bigl(\frac{\textup{d}(\Omega)}{2}\,\left\|\frac{\slo^2}{\sigma}\right\|_{L^\infty(\Omega)}
+\frac{\textup{d}(\Omega)^2}{\gamma}\|\sigma\|_{L^\infty(\Omega)}\Bigr)^2 \|\frac{1}{\slo}\|_{L^\infty(\Omega)}^2  
+ \frac{2}{(1-c_\Omega)} \frac{\textup{d}(\Omega)^2}{\gamma}.$ \\
Together with the first identity in \eqref{enid_Helmholtz}, this yields 
\[
\|\nabla_x u\|_{L^2(\Omega)}^2
+\omega^2\|\slo u\|_{L^2(\Omega)}^2 
\leq \Bigl(2 C +\frac{1}{\omega^2}\left\|\frac{1}{\slo}\right\|_{L^\infty(\Omega)}^2\Bigr)\|f\|_{L^2(\Omega)}^2
\]
and again, elliptic regularity yields \eqref{estHelmholtz2}.
\end{proof}

Analogously to Corollary~\ref{cor:par2state} we obtain well-definedness and differentiability of the parameter-to-state map. 
\begin{corollary}\label{cor:par2state}
Under the assumptions of Corollary~\ref{cor:par2state}, but without \eqref{smallness_perturbedSchr} and with $V$ replaced by 
$V=\Bigl(L^\infty(0,L;H^1(\tilde{\Omega}_\zbot)\cap W^{1,\infty}(0,L;L^2(\tilde{\Omega}_\zbot)\Bigr)^2\times H^2(\Omega)$, the parameter-to-state map
\begin{equation}\label{par2state_secforward}
S:(\slo,\nlc)\mapsto (\phih_1,\phih_2,\psih) \mbox{ solving \eqref{PDEsTransf}, \eqref{BCsTransf}}
\end{equation}
is well-defined as an operator 
$S:W^{1,\infty}(0,L;L^\infty(\tilde{\Omega}_\zbot))\times L^2(\Omega)\to V$
and G\^{a}teaux differentiable for all $(\slo,\nlc)\in (W^{1,\infty}(0,L;L^\infty(\tilde{\Omega}_\zbot))\cap W^{2,q}(0,L;L^\infty(\tilde{\Omega}_\zbot))\times L^2(\Omega)$ as an operator 
$W^{1,\infty}(0,L;L^\infty(\tilde{\Omega}_\zbot))\times L^2(\Omega)\to V$.
\end{corollary}

\section{Convergence of a frozen Newton method for the inverse problem}\label{inverse}
Based on the parameter-to-state map $S$ defined in \eqref{par2state} and analyzed in Section~\ref{sec:forward} as well as the observation operator \eqref{obsop}, we can write the inverse problem as 
\begin{equation}\label{Fqy_inv}
F(\slo,\nlc)=y
\end{equation}
with $y=\pmeas$ and the forward operator $F=C\circ S$ and apply Newton's method to \eqref{Fqy_inv}, using the fact that $F'(\slo,\nlc)=C\circ S'(\slo,\nlc)$ and relying on our verification of G\^{a}teaux differentiability from Section~\ref{sec:forward}.

Alternatively we here consider an all-at-once formulation 
\begin{equation}\label{Fqy_inv_aao}
\mathbb{F}(\slo,\nlc,\phih_1,\phih_2,\psi)=y
\end{equation}
with $y=(0,0,0,0,0,\pmeas)$
based on the variational form of the system of initial-boundary value / boundary value problems \eqref{PDEs}, \eqref{BCs}

\begin{equation}\label{bbF}
\begin{aligned}
\langle\mathbb{F}_k(\vec{q}),\vec{w}\rangle_{\mathbb{X}}:=&
\int_0^L\Bigl\{ \, \int_{\tilde{\Omega}_y} \Bigl(
{\imath\omega_k\slotil_k\,\phih_k \overline{w}_k+}
2\imath\omega_k\slo_k\partial_z\phih_k \,\overline{w}_k
+\nabla_\zbot \phih_k \cdot \nabla_\zbot \overline{w}_k \, \Bigr) d\zbot  \,
+\imath\omega_k\int_{\partial\tilde{\Omega}_\zbot}\sigma_k \phih_k\, \overline{w}_k \, d\Gamma(\zbot)\, \Bigr\} \, dz,\\
\langle\mathbb{F}_3(\vec{q}),\vec{w}\rangle_{\mathbb{X}}:=&
\int_0^L\int_{\tilde{\Omega}_y} \Bigl\{
{\imath\omega_d\slotil_1\,\psih \overline{v}+}
2 \imath \omegad \slo_1 \partial_z\psih\,\overline{v}
+\nabla_\zbot \psih \cdot \nabla_\zbot \overline{v}
+\tilde{\varepsilon}\,\partial_z \psih\, \partial_z \overline{v}
-(\imath\omegap \tilde{\varepsilon}^{-1}\,\nlc \,\phih_1\, \overline{\phih}_2)\,\overline{v}
\Bigr\}\,d\zbot \, dz\\ 
&+\imath\omegad\int_0^L\int_{\partial\tilde{\Omega}_\zbot}\sigma \psih\, \overline{v} \, d\Gamma(\zbot)\, dz
+\imath\omegad\,\tilde{\varepsilon}\,\Bigl(
\sigma_0 \int_{\tilde{\Omega}_y}\psih(0)\, \overline{v}(0) \,d\zbot
+\sigma_L \int_{\tilde{\Omega}_y}\psih(L)\, \overline{v}(L) \,d\zbot\Bigr),\\
\langle\mathbb{F}_{k+3}(\vec{q}),\vec{w}\rangle_{\mathbb{X}}:=&
\int_{\tilde{\Omega}_y}(\phih_k(0) - h_k)\,\overline{u}_k\, d\zbot,
\\
\langle\mathbb{F}_{6}(\vec{q}),\vec{w}\rangle_{\mathbb{X}}:=&
\int_{\Gamma}\imath \omegad\psih \,\overline{u}\, d\Gamma
\end{aligned}
\end{equation}
for $k\in\{1,2\}$, where we have used the abbreviation 
$\vec{q} = (\slo_1,\slo_2,\nlc,\phih_1,\phih_2,\psih)$ and where 
$\vec{w} = (w_1,w_2,v,u_1,u_2,u)$ varies over some space of test functions.
Note that we have introduced a second copy of the variable $\slo$ in order to be able to verify the range invariance condition from \cite{rangeinvar} that allows to prove convergence of a frozen Newton method. 
The additional constraint 
\[
0=\mathcal{P}(\vec{q}):=\slo_1-\slo_2,
\]
that will be imposed via a penalty term during the Newton type iteration \eqref{frozenNewtonHilbert} takes care of merging these two $\slo$ versions into one as the iteration proceeds.

The formal linearization of $\mathbb{F}=(\mathbb{F}_1,\ldots,\mathbb{F}_6)$ at some $\vec{q}^0$ is given by
\begin{equation}\label{bbF_diff}
\begin{aligned}
\langle\mathbb{F}_k'(\vec{q})\dq,\vec{w}\rangle_{\mathbb{X}}=&
\int_0^L\Bigl\{ \int_{\tilde{\Omega}_y} \Big(
\imath\omega_k(
{
\partial_z\dslo_k\phih_k
+\partial_z\slo_k\dphi_k}
+2\dslo_k\partial_z\phih_k 
+2\slo_k\partial_z\dphi_k)
\,\overline{w}_k
+\nabla_\zbot \dphi_k \cdot \nabla_\zbot \overline{w}_k \, \Big) d\zbot\, 
\\ &\hspace*{1cm}
+\imath\omega_k\int_{\partial\tilde{\Omega}_\zbot}\sigma_k \dphi_k\, \overline{w}_k \, d\Gamma(\zbot)\Bigr\} dz,\\
\langle\mathbb{F}_3'(\vec{q})\dq,\vec{w}\rangle_{\mathbb{X}}=&
\int_0^L\int_{\tilde{\Omega}_y} \Bigl\{
\imath \omegad 
({
\partial_z\dslo_1\psih
+\partial_z\slo_1\dpsi}
+2\dslo_1 \partial_z\psih+2\slo_1 \partial_z\dpsi)\,\overline{v}
+\nabla_\zbot \dpsi \cdot \nabla_\zbot \overline{v}
+\tilde{\varepsilon}\,\partial_z \dpsi\, \partial_z \overline{v}\\
&\qquad
-\imath\omegap \tilde{\varepsilon}^{-1}\,
(\dnlc \,\phih_1\, \overline{\phih}_2
+\nlc \,\dphi_1\, \overline{\phih}_2
+\nlc \,\phih_1\, \overline{\dphi}_2
)\,\overline{v}
\Bigr\}\,d\zbot \, dz\\ 
&+\imath\omegad\int_0^L\int_{\partial\tilde{\Omega}_\zbot}\sigma \dpsi\, \overline{v} \, d\Gamma(\zbot)\, dz
+\imath\omegad\,\tilde{\varepsilon}\,\Bigl(
\sigma_0 \int_{\tilde{\Omega}_y}\dpsi(0)\, \overline{v}(0) \,d\zbot
+\sigma_L \int_{\tilde{\Omega}_y}\dpsi(L)\, \overline{v}(L) \,d\zbot\Bigr),\\
\langle\mathbb{F}_{k+3}'(\vec{q})\dq,\vec{w}\rangle_{\mathbb{X}}=&
\int_{\tilde{\Omega}_y}\dphi_k(0)\,\overline{u}_k\, d\zbot,
\\
\langle\mathbb{F}_{6}'(\vec{q})\dq,\vec{w}\rangle_{\mathbb{X}}=&
\imath \omegad\int_{\Gamma}\dpsi\,\overline{u}\, d\Gamma.
\end{aligned}
\end{equation}
It is a bounded linear operator and thus it is the G\^{a}teaux derivative of $\mathbb{F}$ when considered, e.g., as a mapping 
\begin{equation}\label{bbXY}
\begin{aligned}
&\mathbb{F}:\mathbb{X}\to\mathbb{Y}\quad \textup{ with }\\
&\mathbb{X}:={W^{1,p}(0,L;L^p(\tilde{\Omega}_\zbot))\times(L^p(0,L;L^p(\tilde{\Omega}_\zbot)))^2} 
\times (W^{1,\infty}(0,L;L^\infty(\tilde{\Omega}_\zbot))\cap L^2(0,L;H^1(\tilde{\Omega}_\zbot)))^3,\\
&\mathbb{Y}:=(L^2(0,L;H^1(\tilde{\Omega}_\zbot)^*))^2 
\times (H^1((0,L)\times\tilde{\Omega}))^* 
\times (L^2(\tilde{\Omega}))^2 
\times L^2(\Gamma)		
\end{aligned}
\end{equation}
for some $p\in[1,\infty]$.
Note that this slightly differs from the function space setting suggested for the reduced setting by Corollary~\ref{cor:par2state}, but these differences are essential for the verification of convergence conditions for the frozen Newton method, see \eqref{est_r-id} below. Thus, we are here making use of the additional freedom in choosing the function spaces that we gain by using an all-at-once formulation \eqref{Fqy_inv_aao} as compared to a reduced one \eqref{Fqy_inv}.

We can achieve the range invariance relation 
\begin{equation}\label{rangeinvar_diff}
\mathbb{F}(\vec{q})-\mathbb{F}(\vec{q}^0)=\mathbb{F}'(\vec{q}^0)r(\vec{q})
\end{equation}
by setting $r(\vec{q}):=\dq=(\dslo_1,\dslo_2,\dnlc,\dphi_1,\dphi_2,\dpsi)$ with
\[
\begin{aligned}
\dslo_k=& 
{ \frac{\phih_k}{\phih_k^0}\cdot(\slo_k-\slo_k^0)
+ \frac{1}{(\phih_k^0)^2}\int_0^z (\slo_k-\slo_k^0)(z')\Bigl(
\phih_k^0\,\partial_z\phih_k-\phih_k\,\partial_z\phih_k^0\Bigr)(z')\, dz'
}\\
\dnlc=& \frac{\phih_1\, \overline{\phih_2}}{\phih_1^0\, \overline{\phih_2^0}}\cdot(\nlc-\nlc^0)+
\frac{(\phih_1-\phih_1^0)(\overline{\phih_2}-\overline{\phih_2^0})}{\phih_1^0\, \overline{\phih_2^0}}\cdot\nlc^0\\
&{-\frac{\omegad\, \tilde{\varepsilon}}{\omegap \,\phih_1^0\, \overline{\phih_2^0}}
\left(\psih\,\partial_z(\slo_1-\slo_1^0) + 2(\slo_1-\slo_1^0)\, \partial_z\psih
- \psih^0\,\partial_z\dslo_1 - 2 \dslo_1\, \partial_z \psih^0\right)
,}\\
\dphi_k=&\phih_k-\phih_k^0, \qquad \qquad \dpsi=\psih-\psih^0.
\end{aligned}
\]
To this end, we have to choose $\phih_1^0$, $\phih_2^0$, $\psih^0$ such that the denominators in the above expression are bounded away from zero. This is possible, due to the fact that in the all-at-once formulation they need not be PDE solutions corresponding to the coefficients $\dslo_1^0$, $\dslo_2^0$, $\dnlc^0$.
Relying on this, we can also establish closeness of $r$ to the identity in the sense that 
{
\[
r(\vec{q})-(\vec{q}-\vec{q}^0) 
= \bigl(\dslo_1-(\slo_1-\slo_1^0),\dslo_2-(\slo_2-\slo_2^0), \dnlc-(\nlc-\nlc^0), 0, 0, 0\bigr)
\]
where 
\[
\begin{aligned}
&\|\dslo_k-(\slo_k-\slo_k^0)\|_{L^p(L^p)}
\leq\|\tfrac{1}{\phih_k^0}\|_{L^\infty(L^\infty)}
\|\phih_k-\phih_k^0\|_{L^\infty(L^\infty)}
\|\slo_k-\slo_k^0\|_{L^p(L^p)}
\\
&\qquad
+ \|\tfrac{1}{\phih_k^0}\|_{L^\infty(L^\infty)}^2 L^{1/p} \|\slo_k-\slo_k^0\|_{L^p(L^p)}
 \|\phih_k^0\,\partial_z(\phih_k-\phih_k^0)-(\phih_k-\phih_k^0)\,\partial_z\phih_k^0\|_{L^{p*}(L^\infty)}
\end{aligned}
\]
and similarly for $\|\partial_z(\dslo_1-(\slo_1-\slo_1^0)\|_{L^p(L^p)}$, using the identity 
$\partial_z(\dslo_1-(\slo_1-\slo_1^0)=2(\slo_1-\slo_1^0)\,\partial_z(\frac{\phih_k}{\phih^0_k})
+\frac{\phih_k-\phih^0_k}{\phih_k^0}\, \partial_z (\slo_1-\slo_1^0)
+\partial_z(\frac{1}{(\phih_k^0)^2})\, \int_0^z (\slo_k-\slo_k^0)(z')\Bigl(
\phih_k^0\,\partial_z\phih_k-\phih_k\,\partial_z\phih_k^0\Bigr)(z')\, dz'$, 
as well as
\[
\begin{aligned}
&\|\dnlc-(\nlc-\nlc^0)\|_{L^p(L^p)}\leq 
\|\tfrac{1}{\phih_1^0\, \overline{\phih_2^0}}\|_{L^\infty(L^\infty)}\ \cdot
\Bigl\{
\|(\phih_1-\phih_1^0)\overline{\phih_2^0}
+\phih_1(\overline{\phih_2}-\overline{\phih_2^0})\|_{L^\infty(L^\infty)}
\|\nlc-\nlc^0\|_{L^p(L^p)}\\
&\hspace*{4cm}+\frac{\omegad\, \tilde{\varepsilon}}{\omegap}
\Bigl(\|(\psih-\psih^0)\,\partial_z(\slo_1-\slo_1^0)\|_{L^p(L^p)} + 2\|(\slo_1-\slo_1^0)\, (\partial_z\psih-\partial_z\psih^0)\|_{L^p(L^p)}\\
&\hspace*{5cm}+ \|\psih^0 \,\partial_z(\dslo_1-(\slo_1-\slo_1^0)\|_{L^p(L^p)} + 2 \|(\dslo_1-(\slo_1-\slo_1^0))\,\partial_z \psih^0\|_{L^p(L^p)}\Bigr)\Bigr\}
\end{aligned}
\]
hence 
\begin{equation}\label{est_r-id}
\begin{aligned}
\|r(\vec{q})-(\vec{q}-\vec{q}^0)\|_{\mathbb{X}}\leq C \|\vec{q}-\vec{q}^0\|_{\mathbb{X}}^2.
\end{aligned}
\end{equation}
}

Based on the range invariance condition \eqref{rangeinvar_diff} we can rewrite the inverse problem of reconstructing $(\slo,\eta)$ as a combination of an ill-posed linear and a well-posed nonlinear problem
\begin{equation}\label{FP}
\begin{aligned}
&\mathbb{F}'(\vec{q}^0)\hat{r}=y-\mathbb{F}(\vec{q}^0),\\
&r(\vec{q})=\hat{r},\\
&\mathcal{P}\vec{q}=0
\end{aligned}
\end{equation}
for the unknowns $(\hat{r},\vec{q})$. 
Here $\vec{q}^0$ is fixed and $U\subseteq \mathbb{X}$ is a neighborhood of $\vec{q}^0$.

For its regularized iterative solution we consider the frozen Newton method
\begin{equation}\label{frozenNewtonHilbert}
\vec{q}^{n+1} \in \mbox{argmin}_{\vec{q}\in U}
\|\mathbb{F}'(\vec{q}^0)(\vec{q}-\vec{q}^n)+\mathbb{F}(\vec{q}^n)-y^\delta\|_{\mathbb{Y}}^2+\alpha_n\|\vec{q}-\vec{q}^0\|_{\mathbb{X}}^2+\|\mathcal{P}\vec{q}\|_{\mathbb{Z}}^2,
\end{equation}
where $y^\delta\approx y$ is the noisy data, $\alpha_n\to0$ as $n\to\infty$ and $\mathbb{Z}=(L^p(0,L;L^p(\tilde{\Omega}_\zbot)))$.
To work in Hilbert spaces, we set $p=2$.

In order to prove convergence of \eqref{frozenNewtonHilbert}, we additionally need a condition on compatibility of the linear operators $\mathbb{F}'(\vec{q}^0)$ and $\mathcal{P}$, see \cite[Lemma B.1]{rangeinvar}.
Sufficient for this condition is linearized uniqueness with $\slo_1=\slo_2$, that is, triviality of the intersection of the nullspaces 
\begin{equation}\label{NFP0}
\mathcal{N}(\mathbb{F}'(\vec{q}^0)) \, \cap \, \mathcal{N}(\mathcal{P})\ =\{0\}.
\end{equation}
For the model \eqref{PDEs}, \eqref{BCs} described by \eqref{bbF} this would likely require more than one excitation as well as observations at several frequencies, along with the corresponding extension of dependency of $\slo$ in order to allow for range invariance \eqref{rangeinvar_diff}.

\medskip

We here rigorously establish \eqref{NFP0} for the alternative model \eqref{PDEsTransf}, \eqref{BCsTransf} using the Helmholtz equation for the outgoing wave $\psih$ and making the simplifying but still very realistic assumption that propagation of the excitation waves described by $\phih_1$, $\phih_2$ takes place in a homogeneous domain with known and constant sound speed.

The unknown parameters are then the squared slowness $\sslo$ and the nonlinearity coefficient $\nlcc$ and the model reads as 
\[
-\omegad^2 \mathcal{M}[\sslo \psic] + \mathcal{A}\psic +\imath \omegad \mathcal{D}\psic = \imath \omegap \,\nlcc\, f,
\]
where $f$ is given by the excitation waves $\phih_1$, $\phih_2$ that can be precomputed under the assumption of $\slo=\slo_0=\frac{1}{c}$ being known in the first equations of \eqref{PDEsTransf}. 
Moroever, on the boundary $\partial\Omega$ we assume to know $\sslo=\slo_0^2=\frac{1}{c^2}$, and define the operators $\mathcal{A}$, $\mathcal{D}$, $\mathcal{M}$ by  
\begin{equation}\label{eqn:calADM}
\begin{aligned}
&\mathcal{A}u = \Bigl(v\mapsto \int_\Omega\nabla u\cdot\nabla v\, dx +\int_{\partial\Omega} \beta\, u\, v\, ds\Bigr),\\
&\mathcal{D}u = \Bigl(v\mapsto b\int_\Omega\nabla u\cdot\nabla v\, dx +\int_{\partial\Omega} (\sigma+b\beta)u\, v\, ds\Bigr)\,, \\
&\mathcal{M}u = \Bigl(v\mapsto \int_\Omega u\,v\, dx +\int_{\partial\Omega} \frac{\sigma b}{\slo_0^2}\, u\, v\, ds\Bigr).
\end{aligned}
\end{equation}
In our linearized uniqueness proof we will assume that $\mathcal{A}$, 
$\mathcal{D}$ and $\mathcal{M}$ 
are simultaneously diagonalizable.
This holds true, e.g., in the case $\sigma=0$, where $\mathcal{M}=\textup{id}$ and $\mathcal{D}=b\mathcal{A}$.
Note that in the all-at-once setting considered here, we are not bound to well-posedness theory of the underlying PDE problems and therefore have more freedom in  choosing the coefficients in the boundary conditions. 

To obtain uniqueness, we will need observations on an interval $I$ of difference frequencies $\omegad$ and with at least two different sets of excitation frequencies, more precisely 
\begin{equation}\label{omegad_kappa}
\begin{aligned}
&\omega_2=\frac{\omega_1\kappa}{\omega_1+\kappa} \, \kappa, \quad\kappa\in\{\kappa_1,\kappa_2\}\subseteq\mathbb{R}^+, \quad \omegad=\omega_1-\omega_2\in I\subseteq\mathbb{R}^+\\
& \textup{with $I$ countable containing at least one accumumation point} 
\end{aligned}
\end{equation}
which implies $\omegap=\omegad^2 \kappa$.
Correspondingly, we consider $\psic=\psic(\omegad,\kappa)$ as a function of $\omegad$, while the original coefficients $\sslo$, $\nlcc$ of course stay independent of $(\omegad,\kappa)$. In order to satisfy range invariance, we will introduce an artificial dependence $\sslo=\sslo(\omegad,\kappa)$, while keeping $\nlcc$ independent of $\omegad$ and $\kappa$.
 
Therewith, the inverse problem reads as 
\[ 
\mathbb{F}(\sslo,\check{\nlc},\psic )=y
\]
with 
\[
\begin{aligned}
&\mathbb{F}_1(\vec{q}) = 
-\omegad^2\,\bigl(\mathcal{M}(\sslo\, \psic)  + \imath \kappa \, \nlcc\, f\bigr) +\mathcal{A} \psic +\imath \omegad \mathcal{D}\psic ,
\\
&\mathbb{F}_2(\vec{q}) = \imath \omegad \textup{tr}_{\check{\Gamma}}\psic,
\end{aligned}
\]
where $\vec{q} = (\sslo,\nlcc,\psic)$, 
and $y=(0,\pmeas)$.

It is readily checked that range invariance \eqref{rangeinvar_diff} holds with 
\[
r(\sslo,\nlcc,\psic) = 
\Bigl((\sslo-\sslo_0)\frac{\psic}{\psic_0},\nlcc-\nlcc_0,\psic-\psic_0\Bigr)
\]
and 
\begin{equation}\label{rid}
\begin{aligned}
\|r(\vec{q})-(\vec{q}-\vec{q}^0)\|_{\mathbb{X}}
\leq \left\|\frac{1}{\psic_0}\right\|_{L^\infty(\Omega)} \hspace*{-0.5cm}\cdot\|(\sslo-\sslo_0)(\psic-\psic_0)\|_{L^2(\Omega)}
\leq C \|\vec{q}-\vec{q}^0\|_{\mathbb{X}}^2 \leq c \|\vec{q}-\vec{q}^0\|_{\mathbb{X}}
\end{aligned}
\end{equation}
with $\mathbb{X}=L^2(\Omega)\times L^2(\Omega)\times H^2(\Omega)$, provided $\|\frac{1}{\psic_0}\|_{L^\infty(\Omega)}<\infty$, that is, $\psic_0$ is bounded away from zero.
Here $c$ can be made small by restricting $\vec{q}$ to a sufficiently small neighborhood of $\vec{q}^0$.

Defining $\mathcal{P}$, $\mathbb{X}$, $\mathbb{Y}$, $\mathbb{Z}$ by 
\begin{equation}\label{PXYZ}
\begin{aligned}
&\mathcal{P}(\vec{q})(\omegad):=\Bigl(\sslo(\omegad,\kappa_1)-\sslo(\omegad,\kappa_2),
\ \sslo(\omegad,\kappa_1)-\frac{1}{\mu(I)}\int_I \sslo(\cdot,\kappa_1)\,d\mu\Bigr),\\
&\mathbb{X}=L^2_\mu(I;L^2(\Omega))^2\times L^2(\Omega)\times H^2(\Omega), \quad
\mathbb{Y}=L^2_\mu(I;L^2(\Gamma))^2, \quad
\mathbb{Z}=L^2_\mu(I;L^2(\Omega))^2
\end{aligned}
\end{equation}
for some finite measure $\mu$ on $I$,
we can write the inverse problem as \eqref{FP} and use \eqref{frozenNewtonHilbert} for its regularized numerical solution.

To prove convergence of \eqref{frozenNewtonHilbert} we require linearized uniqueness \eqref{NFP0}, which we verify as follows. For $(\dsslo,\dnlcc,\dpsic)\in\mathcal{N}(\mathcal{P})$, that is, $\sslo$ independent of $\omegad$ and $\kappa$, the assumption $\mathbb{F}'(\sslo_0,\nlcc_0,\psic_0)(\dsslo,\dnlcc,\dpsic)=0$ with $\sslo_0:\equiv\frac{1}{c^2}$ and $\nlcc_0$ arbitrary reads as 
\[
\begin{aligned}
(-\omegad^2\,\tfrac{1}{c^2}\mathcal{M}+ \mathcal{A} +\imath \omegad \mathcal{D})\, \dpsic(\omegad,\kappa)
&=\omegad^2 \bigl(\mathcal{M}[\dsslo\, \psic_0(\omegad,\kappa)]  + \imath\kappa \, \dnlcc\, f(\omegad,\kappa)\bigr),\\
\imath \omegad \textup{tr}_{\check{\Gamma}}\dpsic(\omegad,\kappa)&=0\\
&\hspace*{2.5cm}\textup{for all }\omegad\in I, \ \kappa\in\{\kappa_1,\kappa_2\}.
\end{aligned}
\]
We now assume that $\psic_0$ and $f$ have been chosen independent of $\omegad$ and $\kappa$ 
\begin{equation}\label{fpsic0}
\psic_0(\omegad,\kappa)\equiv\psic_0, \quad f(\omegad,\kappa)\equiv f
\textup{ with }\psic\not=0, \ f\not=0\textup{ a.e.}
\end{equation}
and diagonalize the operators $\mathcal{A}$, 
$\mathcal{D}$, $\mathcal{M}$, 
by means of their eigensystems $(\lambda_j,\varphi_j^k)_{j\in\mathbb{N},k\in K^j}$, 
$(\rho_j,\varphi_j^k)_{j\in\mathbb{N},k\in K^j}$ 
$(\mu_j,\varphi_j^k)_{j\in\mathbb{N},k\in K^j}$ 
--- recall that we have assumed them to be jointly diagonalizable and that they are symmetric and positive semidefinite.
Note that with $\slo_0>0$, $\sigma,\,\beta,\,b\geq0$, we have 
\[
\mu_j\geq 1, \quad \lambda_j,\, \rho_j\geq 0, \quad \lambda_j\to\infty \textup{ as }j\to\infty. 
\]
We will additionally assume
\begin{equation}\label{eqn:rho_lambda_mu}
\Bigl(\frac{\rho_j}{\mu_j}=\frac{\rho_\ell}{\mu_\ell} \textup{ and } 
\frac{\lambda_j}{\mu_j^2}=\frac{\lambda_\ell}{\mu_\ell^2}\Bigr)
\ \Rightarrow \ j=\ell,
\end{equation} 
which is the case, e.g., if $\sigma=0$, since then $\mu_j\equiv1$.

This yields
\begin{equation}\label{ajkbjk0}
\begin{aligned}
0=\dpsic(x_0;\omegad,\kappa)=c^2\omegad^2\sum_{j=1}^\infty\sum_{k\in K^j}
\frac{1}{-\omegad^2\,\mu_j+c^2\,\lambda_j+\imath\omegad c^2 \,\rho_j} (a_j^k + \imath\kappa b_j^k)\,
\varphi_j^k(x_0)\\
\textup{ for all }x_0\in\Gamma,\ \omegad\in I, \ \kappa\in\{\kappa_1,\kappa_2\}
\end{aligned}
\end{equation}
for 
\[
a_j^k = \mu_j\langle \dsslo\, \psic_0,\varphi_j^k\rangle_{L^2(\Omega)}, \quad
b_j^k = \langle \dnlcc\, f,\varphi_j^k\rangle_{L^2(\Omega)}.
\]
It is straightforward to show that the functions $\omegad\mapsto\frac{1}{w_j(\imath\omegad)}$ with $w_j(\imath\omegad):=-\omegad^2\,\mu_j+c^2\,\lambda_j+\imath\omegad c^2 \,\rho_j$ are linearly independent on $I$. 
Indeed, assuming $\sum_{j=1}^\infty \frac{1}{w_j(\imath\omegad)}\, c_j=0$ for all $\omegad\in I$ and a sequence $\vec{c}=(c_j)_{j\in\mathbb{N}}$, after multiplication with $\prod_{\ell\in\mathbb{N}}w_\ell(\imath\omegad)$, we conclude  $0=W^{\vec{c}}(z):=\sum_{j=1}^\infty\prod_{\ell\not=j}w_\ell(z)\, c_j$ for all $z\in\imath I$ and thus, since the function $W^{\vec{c}}$ is analytic, it vanishes on all of $\mathbb{C}$. Inserting the roots 
$z^k_\pm=-\frac{c\rho_j}{2}\pm\sqrt{\frac{c^2\rho_j^2}{4}-\lambda_j}$ of $w_k$ and using the fact that under condition \eqref{eqn:rho_lambda_mu} they are distinct (in the sense that $z^j_+=z^\ell_+$ and $z^j_-=z^\ell_-$ implies $j=\ell$), we obtain $\vec{c}=0$. 
\\
By the linear independence of $(\frac{1}{w_j})_{j\in\mathbb{N}}$ we conclude from \eqref{ajkbjk0} 
\[
\sum_{k\in K^j} (a_j^k + \imath\kappa b_j^k)\, \varphi_j^k(x_0) = 0 \quad
\textup{ for all }x_0\in\Gamma,\ j\in\mathbb{N}, \ \kappa\in\{\kappa_1,\kappa_2\}.
\]
Under a linear independence assumption on the individual eigenspaces (cf., e.g., \cite{nonlinearityimaging_fracWest})  
\begin{equation}\label{eqn:ass_inj_Gamma_rem}
\textup{ for all } j\in\mathbb{N}\ : \quad 
\left(\sum_{k\in K^j} \gamma_k \varphi_j^k(x) = 0 \ \mbox{ for all }x\in\Gamma\right)
\ \Longrightarrow \ \left(\gamma_k=0 \mbox{ for all }k\in K^j\right)
\end{equation}
we conclude
\[
\begin{cases}
a_j^k + \imath\kappa_1 b_j^k=0\\
a_j^k + \imath\kappa_2 b_j^k=0
\end{cases}
\quad \textup{ for all }j\in\mathbb{N}, \ k\in K^j, 
\]
which implies $a_{j}^k=0$, $b_j^k=0$, for all $j\in\mathbb{N}$, $k\in K^j$, thus under assumption \eqref{fpsic0}, we get $\dsslo=0$, $\dnlcc=0$.

According to \cite[Theorem 2]{rangeinvar}, we obtain the following.
\begin{theorem}\label{thm:convfrozenNewton}
Let  $\vec{q}^\dagger=(\sslo,\nlcc,\psic)$ be a solution to \eqref{FP} and let for the noise level $\delta\geq\|y^\delta-y\|_{\mathbb{Y}}$ the stopping index $n_*=n_*(\delta)$ be chosen such that 
\begin{equation}\label{nstar}
n_*(\delta)\to0, \quad \delta\sum_{j=0}^{n_*(\delta)-1}c^j\alpha_{n_*(\delta)-j-1}^{-1/2} \to 0 \qquad \textup{ as }\delta\to0
\end{equation}
with $c$ as in \eqref{rid}.
Moreover, let \eqref{fpsic0} with $\frac{1}{\psic_0}\in L^\infty(\Omega)$ hold and let the operators 
$\mathcal{A}$, $\mathcal{D}$, $\mathcal{M}$ be jointly diagonalizable with eigenvalues satisfying \eqref{eqn:rho_lambda_mu}.

Then there exists $\varrho>0$ sufficiently small such that for $\vec{q}_0\in\mathcal{B}_\varrho(\vec{q}^\dagger)$ the iterates $(\vec{q}_n^\delta)_{n\in\{1,\ldots,n_*(\delta)\}}$ are well-defined by \eqref{frozenNewtonHilbert}, 
remain in $\mathcal{B}_\varrho(\vec{q}^\dagger)$ and converge in $\mathbb{X}$, $\|\vec{q}_{n_*(\delta)}^\delta-\vec{q}^\dagger\|_{\mathbb{X}}\to0$ as the noise level $\delta\to0$. In the noise free case $\delta=0$, $n_*(\delta)=\infty$ we have $\|\vec{q}_n-\vec{q}^\dagger\|_{\mathbb{X}}\to0$ as $n\to\infty$.
\end{theorem}



\begin{funding}
This work was supported by the Austrian Science Fund {\sc fwf}  http://dx.doi.org/10.13039/501100002428 under the grant {\sc doc}78.
\end{funding}

\end{document}